\documentclass[10pt, reqno]{article}
\usepackage{latexsym, amsmath, amssymb, a4, epsfig}

\newcommand{\eqnref}[1]{(\ref {#1})}
\renewcommand{\ker}{\operatorname{Ker}}
\newcommand{\coker}{\operatorname{Coker}}
\newcommand{\ind}{\operatorname{Ind}}
\newcommand{\tr}{\operatorname{tr}}
\newcommand{\im}{\operatorname{Im}}
\newcommand{\rank}{\operatorname{rank}}
\newcommand{\rot}{\operatorname{rot}}

\newcommand{\ds}{\displaystyle}

\newcommand{\p}{\partial}
\newcommand{\pf}{\medskip \noindent {\sl Proof}. \ }
\newcommand{\qed}{\hfill $\Box$ \medskip}
\newcommand{\RR}{\mathbb{R}}
\newcommand{\ZZ}{\mathbb{Z}}
\newcommand{\Scal}{\mathcal{S}}
\newcommand{\Kcal}{\mathcal{K}}
\newcommand{\Acal}{\mathcal{A}}
\newcommand{\Dcal}{\mathcal{D}}
\newcommand{\Qcal}{\mathcal{Q}}

\newcommand{\tnu}{\widetilde{\nu}}
\newcommand{\tmu}{\widetilde{\mu}}
\newcommand{\tlambda}{\widetilde{\lambda}}

\newcommand{\vp}{\varphi}
\def\nm{\noalign{\medskip}}
\newcommand{\pd}[2]{\frac {\p #1}{\p #2}}
\newcommand{\beq}{\begin{equation}}
\newcommand{\eeq}{\end{equation}}

\newcommand{\Om}{\Omega}

\newtheorem{thm}{Theorem}[section]
\newtheorem{cor}[thm]{Corollary}
\newtheorem{lem}[thm]{Lemma}

\numberwithin{equation}{section}

\begin{document}
\title{Asymptotic Analysis of High-Contrast Phononic Crystals and a Criterion for the Band-Gap Opening\thanks{H.A. is partially
supported by the Brain Pool Korea Program at Seoul National
University, H.K. is partially supported by the grant KOSEF
R01-2006-000-10002-0, and H.L. is supported by BK21 Math. division
at Seoul National University.}}
\author{Habib Ammari\thanks{Centre
de Math\'ematiques Appliqu\'ees, CNRS UMR 7641 and Ecole
Polytechnique,
    91128 Palaiseau Cedex, France
   (habib.ammari@polytechnique.fr).} \and Hyeonbae Kang\thanks{School of Mathematical
Sciences, Seoul National University, Seoul 151-747, Korea
(hkang@math.snu.ac.kr, hdlee@math.snu.ac.kr).}  \and Hyundae
Lee\footnotemark[3]} \maketitle

\begin{abstract}

We investigate the band-gap structure of the frequency spectrum
for elastic waves in a high-contrast, two-component periodic
elastic medium. We consider two-dimensional phononic crystals
consisting of a background medium which is perforated by an array
of holes periodic along each of the two orthogonal coordinate
axes. In this paper we establish a full asymptotic formula for
dispersion relations of phononic band structures as the contrast
of the shear modulus and that of the density become large. The
main ingredients are integral equation formulations of the
solutions to the harmonic oscillatory linear elastic equation and
several theorems concerning the characteristic values of
meromorphic operator-valued functions in the complex plane such as
Generalized Rouch\'{e}'s theorem. We establish a connection
between the band structures and the Dirichlet eigenvalue problem
on the elementary hole. We also provide a criterion for exhibiting
gaps in the band structure which shows that smaller the density of
the matrix is, wider the band-gap is, provided that the criterion
is fulfilled. This phenomenon was reported by Economou and Sigalas
in \cite{ES94} who observed that periodic elastic composites whose
matrix has lower density and higher shear modulus compared to
those of inclusions yield better open gaps. Our analysis in this
paper agrees with this experimental finding.

\vskip 0.5\baselineskip
\end{abstract}

\tableofcontents

\section{Introduction}
In the past decade there has been a steady growth of interest in
the motion of elastic waves through inhomogeneous materials. The
primary motive for these investigations has been the design of the
so-called phononic band gap materials or phononic crystals. The
most recent research in this field has focused on theoretical and
experimental demonstration of band gaps in two-dimensional and
three-dimensional structures constructed of high-contrast elastic
materials arranged in a periodic array. This type of structure
prevents elastic waves in certain frequency ranges from
propagating and could be used to generate frequency filters with
control of pass or stop bands, as beam splitters, as sound or
vibration protection devices, or as elastic waveguides. See, for
example, \cite{vasseur, cheng2, Kushwa, Sigalas}.

The interest in phononic crystals has been renewed by the work in
optics on photonic band gap materials. Since the seminal paper by
Yablonovitch \cite{yabl}, significant progress has been made in
microstructuring a dielectric or magnetic material on the scale of
the optical wavelength so that a range of frequencies for which
incident electromagnetic waves are unable to propagate through the
designed crystal exists.  See \cite{joa2, FigoKuch2, FigoKuch3, Kuch2}. See
also \cite{yuri} for extensive list of references on photonic
crystals.

To formulate the investigation of this paper in mathematical
terms, let $D$ be a connected domain with the Lipschitz boundary
lying inside the open square $]0, 1[^2$. An important example of
phononic crystals consists of a background elastic medium of
constant Lam\'e parameters $\lambda$ and $\mu$ which is perforated
by an array of arbitrary-shaped inclusions $\Omega = \cup_{n\in
\mathbb{Z}^2}(D+n)$ periodic along each of the two orthogonal
coordinate axes in the plane. These inclusions have Lam\'{e}
constants $\tilde{\lambda}$, $\tilde{\mu}$. The shear modulus
$\mu$ of the background medium is assumed to be larger than that
of the inclusion $\tilde\mu$. Then we investigate the spectrum of
the self-adjoint operator defined by
\begin{equation} \label{periodicop}
\mathbf{u} \mapsto - \nabla\cdot(C\nabla\mathbf{u}) = -
\sum_{j,k,l=1}^2 \frac{\partial}{\partial x_j} \bigg( C_{ijkl}
\frac{\partial u_k}{\partial x_l} \bigg),
\end{equation}
which is densely defined on $L^2(\RR^2)^2$. Here  the elasticity
tensor $C$ is given by \beq C_{ijkl} := \ds \bigg ( \lambda \,
\chi(\RR^2 \setminus \overline \Omega) + \tlambda \, \chi (\Omega)
\bigg) \delta_{ij} \delta_{kl} + \bigg( \mu \, \chi(\RR^2
\setminus \overline \Omega) + \tmu\, \chi (\Omega) \bigg)
(\delta_{ik} \delta_{jl} + \delta_{il} \delta_{jk} ), \eeq where
$\chi(\Om)$ is the indicator function of $\Om$.

In this paper we adopt this specific two-dimensional model to
understand the relationship between the contrast of the shear modulus and the band gap structure of the phononic crystal.
We will also consider the case of two materials with different densities in order to investigate the relation
between the density contrast and the band gap.

By Floquet theory \cite{KUCH}, the spectrum of the Lam\'{e} system with
periodic coefficients is represented as a union of bands, called
phononic band structure. Carrying out a band structure calculation
for a given phononic crystal involves a family of eigenvalue
problems, as the quasi-momentum is varied over the first Brillouin
zone. The problem of finding the spectrum of \eqnref{periodicop}
is reduced to a family of eigenvalue problems with
quasi-periodicity condition, {\it i.e.},
\begin{align}\label{eigenperiodic}
\nabla\cdot(C\nabla\mathbf{u})+\omega^2
\mathbf{u}=0~\mbox{in}~\RR^2,
\end{align}
with the periodicity condition
 \beq
 \mathbf{u}(x+n)=e^{i\alpha\cdot n}\mathbf{u}(x) \quad \mbox{for every }
 n\in\ZZ^2.
 \eeq
Here the quasi-momentum $\alpha$ varies over the Brillouin zone
$[0,2\pi[^2$. Each of these operators has compact resolvent so
that its spectrum consists of discrete eigenvalues of finite
multiplicity. We show that these eigenvalues are the
characteristic values of meromorphic operator-valued functions
that are of Fredholm type with index zero. This yields a new and
natural approach to the computation of the band gap phononic
structure which is based on a combination of boundary element
methods and Muller's method \cite{muller} for finding complex
roots of scalar equations. Efficiency of a similar approach for
computing photonic band gaps has been demonstrated in
\cite{cheng1,cheng2}. We then proceed from the generalized
Rouch{\'e}'s theorem to construct their complete asymptotic
expressions as the Lam\'e parameter $\mu$ of the background goes
to infinity. For $\alpha\ne 0$, we prove that the discrete
spectrum of \eqnref{eigenperiodic} accumulates near the Dirichlet
eigenvalues of Lam\'e system in $D$ as $\mu$ goes to infinity. We then
obtain a full asymptotic formula for the eigenvalues with the
leading order term of order $\mu^{-1}$ calculated explicitly.  For
the periodic case $\alpha=0$, we establish a formula for asymptotic behavior
of eigenvalues, but their limiting set is generically different
from that for $\alpha\ne 0$. We also consider the case when
$|\alpha|$ is of order $1/ \sqrt{\mu}$ and derive an asymptotic
expansion for the eigenvalues in this case as well. Not surprisingly, this
formula tends continuously to the previous ones as $\alpha
\sqrt{\mu}$ goes to zero or to infinity. We finally provide a
criterion for exhibiting gaps in the band structure. As has been
said, the existence of those spectral gaps implies that the
elastic waves in those frequency ranges are prohibited from
travelling through the elastic body. Our criterion shows that
smaller the density of the matrix is, wider the band-gap is,
provided that the criterion is fulfilled. This phenomenon was reported by
Economou and Sigalas in \cite{ES94} who observed that periodic
elastic composites whose matrix has lower density and higher shear
modulus compared to those of inclusions yield better open gaps.

Similar results for the photonic crystals were obtained by Hempel
and Lienau in \cite{hempel2}, where they dealt with conductivity
equation with high contrast in two phase composites. See also
Friedlander \cite{friedlander2}. Another related work is
\cite{AKSZ} which concerns the photonic band gap using the same
method as the one in this work. A justification of the existence
of elastic band gaps in periodic composite materials with strong
heterogeneities has been recently provided by \'Avila et al. in
\cite{avila} by extending Bouchitt\'e and Felbacq scalar
homogenization approach \cite{bouchitte} to the elasticity problem. We also
mention works by Movchan and his collaborators \cite{mcphedran,
poul, plats, mm}. To the best of our knowledge, our result on the gap
opening is achieved by using a method significantly different from
those in the literature. All the other asymptotic results of this paper are new
and have never been established elsewhere.

The main ingredients in deriving the results of this paper are the
boundary integral equations and the theory of meromorphic
operator-valued functions. Using integral representations of the
solutions to the harmonic oscillatory linear elastic equation, we
reduce this problem to the study of characteristic values of
integral operators in the complex planes. Generalized Rouch\'{e}'s
theorem and other techniques from the theory of meromorphic
operator-valued functions are combined with careful asymptotic
expansions of integral kernels to obtain full asymptotic
expansions for eigenvalues. This method was first used in
\cite{AT}, and then successfully applied to obtain an asymptotic
formula for the eigenvalues of Laplacian under singular
perturbations \cite{AKLZ} and high contrast asymptotics for the
photonic crystals \cite{AKSZ}. See also \cite{AKL2}.

Results of this paper could be used to design an optimization tool
based on layer potential techniques for the systematic design of
the band-gap elastic materials and structures. Since the limiting
situation reduces to easy-to-calculate spectra, the idea would be
to start with these spectra (as initial guess) and then compute
the gradient of some target functional using our asymptotic
expansions with respect to the contrast. Moreover, in order to
optimize the position and width of these gaps, we only need to
optimize the shape of the inclusion considering the (more simpler)
limiting situation.

The paper is organized as follows. In section~2, we recall
relevant definitions and state the generalized Rouch\'{e}'s
theorem. Then we introduce the single and the double layer
potentials for the harmonic oscillatory linear elastic equation
and cover well-known results for the quasi-periodic case as well.
In section~3, we state and prove our main results. Section~4 is
devoted to the derivation of a criterion for gap opening in the
spectrum of the operator given by \eqnref{periodicop} as $\mu
\rightarrow + \infty$. In section 5, we derive a similar criterion when the contrast of the density is high.

\section{Preliminaries}

The first subsection of this section covers relevant parts of the theory of
meromorphic operator-valued functions. We state the generalized
Rouch\'e's theorem. In the second subsection we review some
well-known results on the solvability and theory of layer potentials for the
harmonic oscillatory linear elastic equation, which we shall use
frequently throughout this paper. The third subsection is devoted
to the radiation condition for elastic wave propagation. In the
fourth subsection we collect some notation and well-known results
regarding quasi-periodic layer potentials for the harmonic
oscillatory linear elastic equation.

\subsection{The generalized Rouch\'{e}'s theorem}
In this work the approach we develop is a boundary integral
technique with rigorous justification based on the generalized
Rouch{\'e}'s theorem. For readers' convenience we recall this the
main results of Gohberg and Sigal in \cite{GS}. We begin by
collecting some definitions.

 Let $\mathcal{G}$ and $\mathcal{H}$ be two Banach spaces and let
$\mathcal{L}(\mathcal{G},\mathcal{H})$ be the set of all bounded
operators from $\mathcal{G}$ to $\mathcal{H}$.
$A\in\mathcal{L}(\mathcal{G},\mathcal{H})$ is called a {\it
Fredholm operator} if $A$ is closed and $\ker A$ and $\coker A$
are finite dimensional. If $A$ is Fredholm, the index of $A$ is
defined to be
\begin{align}
\ind A = \dim \ker A - \dim \coker A.
\end{align}
If $A$ is Fredholm, then for any compact operator
$K\in\mathcal{L}(\mathcal{G},\mathcal{H})$, $A+K$ is also Fredholm
with index $\ind A$.

Let $U$ be an open set in $\mathbb{C}$. Suppose that
$\omega\mapsto A(\omega)$ is an operator-valued function from $U$
into $\mathcal{L}(\mathcal{G},\mathcal{H})$. We call $\omega_0\in
U$ a {\it characteristic value} of $A(\omega)$ if
\begin{itemize}
\item $A(\omega)$ is holomorphic in some neighborhood of
$\omega_0$, except possibly for $\omega_0$; \item There exists a
function $\phi(\omega)$ from a neighborhood of $\omega_0$ to
$\mathcal{G}$, holomorphic and nonzero at $\omega_0$, such that
$A(\omega)\phi(\omega)$ is holomorphic at $\omega_0$ and
$A(\omega_0)\phi(\omega_0)=0$.
\end{itemize}
The function $\phi(\omega)$ in the above definition is called a
{\it root function} of $A(\omega)$ associated to $\omega_0$ and
$\phi(\omega_0)$ is called an {\it eigenvector}. The closure of
the space of eigenvectors corresponding to $\omega_0$ is denoted
by $\ker A(\omega_0)$.

Let $\phi_0$ be an eigenvector corresponding to $\omega_0$. Let
$V(\omega_0)$ be a complex neighborhood of $\omega_0$. The {\it
rank} of $\phi_0$ is the largest integer $m$ such that there exist
$\phi : V(\omega_0)\rightarrow\mathcal{G}$ and $\psi :
V(\omega_0)\rightarrow\mathcal{H}$ holomorphic satisfying
\begin{align}
A(\omega)\phi(\omega)=(\omega-\omega_0)^m \psi(\omega), \quad \phi(\omega_0)=\phi_0, ~~\psi(\omega_0)\ne 0.
\end{align}

Suppose that  $A(\omega)$ admits the Laurent expansion in
$\omega_0$ such as
\begin{align}
A(\omega)=\sum_{j\geq -s}(\omega-\omega_0)^j A_j
\end{align}
for some non-negative integer $s$, where $A_j$ are operators in
$\mathcal{L}(\mathcal{G},\mathcal{H})$ independent of $\omega$. If
the operators $A_j$, $j=-s,\cdots,-1$, are finite dimensional,
then $A(\omega)$ is called {\it finitely meromorphic} at
$\omega_0$. If the operator $A_0$ is a Fredholm operator,
$A(\omega)$ is said to be of {\it Fredholm type} at $\omega_0$.

If $A(\omega)$ is holomorphic at $\omega_0$ and $A(\omega_0)$ is
invertible, then $\omega_0$ is called a {\it regular point } of
$A(\omega)$. A point $\omega_0$ is called a {\it normal point} of
$A(\omega)$ if $A(\omega)$ is finitely meromorphic and of Fredholm
type at $\omega_0$ and there exists some neighborhood
$V(\omega_0)$ of $\omega_0$ in which all the points except
$\omega_0$ are regular points of $A(\omega)$. The following holds.
\begin{lem}\label{normalpt}
Every normal point $\omega_0$ of $A(\omega)$ is a normal point of
$A^{-1}(\omega)$.
\end{lem}

Suppose that $\omega_0$ is a normal point of $A(\omega)$. Then
$\dim \ker A(\omega_0)< + \infty$ and the ranks of all vectors in
$\ker A(\omega_0)$ are finite. A system of eigenvectors $\phi^j$,
$j=1,\cdots,n$, is called a {\it canonical system of eigenvectors}
of $A(\omega)$ associated to $\omega_0$ if
$\{\phi^1,\ldots,\phi^n\}$ is a basis of $\ker A(\omega_0)$ and
the rank of $\phi^j$ is the maximum of the ranks of all
eigenvectors in some direct complement in $\ker A(\omega_0)$ of
the linear span of the vectors $\phi^1, \cdots , \phi^{j-1}$. Then
we define the {\it null multiplicity }of the characteristic value
$\omega_0$ of $A(\omega)$ by
\begin{align}
N(A(\omega_0))=\sum_{i=1}^n \rank (\phi^i).
\end{align}
Similarly, we can also define $N(A^{-1}(\omega_0))$. Then the
number
\begin{align}
M(A(\omega_0))=N(A(\omega_0))-N(A^{-1}(\omega_0))
\end{align}
is called the {\it multiplicity} of the characteristic value
$\omega_0$ of $A(\omega)$.

For an open subset $V$ of $U$, an operator-valued function
$A(\omega)$ is called {\it normal} with respect to $\p V$ if
$A(\omega)$ is invertible and continuous at $\p V$ and all the
points of $V$ are regular, except for a finite number of  normal
points of $A(\omega)$.

Suppose that $A(\omega)$ is normal with respect to $\p V$ and let
$\omega_i$, $i=1,\cdots, \sigma$, be a finite number of normal
points of $A(\omega)$ in $V$. Then we define
\begin{align}
\mathcal{M}(A(\omega); \p V)=\sum^\sigma_{i=1}M(A(\omega_i)).
\end{align}

The generalization
of Rouch\'{e}'s theorem was obtained in \cite{GS}.

\begin{thm}\label{rouche}
Let $A(\omega)$ be an operator-valued function which is normal
with respect to $\p V$. If an operator-valued function $S(\omega)$
which is finitely meromorphic in $V$ and continuous on $ \p V$
satisfies the condition
\begin{align}
\|A^{-1}(\omega)S(\omega)\|_{\mathcal{L}(\mathcal{G},\mathcal{G})}<1,
~~ \omega\in \p V,
\end{align}
then $A(\omega)+S(\omega)$ is also normal with respect to $\p V$ and
\begin{align}
\mathcal{M}(A(\omega); \p V)=\mathcal{M}(A(\omega)+S(\omega); \p V).
\end{align}
\end{thm}

The generalization of Steinberg theorem is given by the following.
\begin{thm}
Suppose that $A(\omega)$ is an operator-valued function which is
finitely meromorphic and of Fredholm type in $ V$. If $A(\omega)$
is invertible at one point of $V$, then $A(\omega)$ has a bounded
inverse for all $\omega\in V$, except possibly for certain
isolated points.
\end{thm}

The operator generalization of the residue theorem is as follows.
\begin{thm}\label{formula}
Suppose that $A(\omega)$ is an operator-valued function which is
normal with respect to $\p V$. Let $f(\omega)$ be a scalar
function which is holomorphic in $V$ and continuous on $\p V$.
Then we have
\begin{align}
\frac{1}{2\pi \sqrt{-1}}\tr\int_{\p V}
f(\omega)A^{-1}(\omega)\frac{d}{d\omega}A(\omega)
d\omega=\sum^\sigma_{j=1}M(A(\omega_j))f(\omega_j),
\end{align}
where $\omega_j$, $j=1,\cdots,\sigma$, are all the points in $ V$
which are either poles or characteristic values of $A(\omega)$.
\end{thm}

Here tr denotes the trace of operator which is the sum of all its
nonzero eigenvalues. The following property of the trace is of use
to us:
\begin{align}\label{trace}
\tr\int_{\p V(\omega_0)}A(\omega)B(\omega)d\omega=\tr\int_{\p
V(\omega_0)}B(\omega)A(\omega)d\omega,
\end{align}
where $A(\omega)$ and $B(\omega)$ are operator-valued functions
which are finitely meromorphic in $V(\omega_0)$, which contains no
poles of $A(\omega)$ and $B(\omega)$ other than $\omega_0$.

\subsection{Green's tensor and layer potentials}
 Let $\Om$ be a bounded domain in $\mathbb{R}^2$ with a
connected Lipschitz boundary $\partial \Omega$. Let the space $H^1
(\partial \Omega)$ be the set of functions $f \in L^2(\partial
\Omega)$ such that ${\partial f}/{\partial T} \in L^2(\partial
\Omega),$ where $\partial / \partial T$ denotes the tangential
derivative on $\partial \Omega$. We use $H^s(\Omega)$ as a
notation for the standard Sobolev space of order $s$.

Let $\lambda$, $\mu$ be the Lam\'{e} constants for $\Om$
satisfying
\begin{equation}
\mu>0~~\mbox{and}~~3 \lambda+ 2 \mu>0.
\end{equation}
The elastostatic system corresponding to the Lam{\'e} constants
$\lambda, \mu$ is defined by
\begin{equation} \label{lamedef}
\mathcal{L}^{\lambda,\mu}\mathbf{u}=\mu\triangle
\mathbf{u}+(\lambda+\mu)\nabla\nabla\cdot \mathbf{u},
\end{equation}
and the corresponding conormal derivative ${\partial
\mathbf{u}}/{\partial \nu}$ on $\partial \Omega$ is defined to be
\begin{equation} \label{lamedefnu}
\pd{\mathbf{u}}{\nu}=\lambda(\nabla\cdot \mathbf{u})N+\mu(\nabla
\mathbf{u} +\nabla \mathbf{u}^t)N,
\end{equation}
where $N$ is the outward unit normal to $\partial \Omega$ and the superscript
$t$ denotes the transpose of a matrix.

The fundamental solution $\mathbf{\Gamma}^\omega =
({\Gamma}^\omega_{ij})_{i,j=1}^2$ to the operator
$\mathcal{L}^{\lambda,\mu}+\omega^2$, $\omega\ne0$, is given by
\begin{equation}\label{kupradze}
{\Gamma}^\omega_{ij}(x)=
\ds-\frac{\sqrt{-1}}{4\mu}\delta_{ij}H^{(1)}_0\left(\frac{\omega|x|}{c_T}\right)
+\frac{\sqrt{-1}}{4\omega^2}\partial_i\partial_j\left(H^{(1)}_0\left(\frac{\omega|x|}{c_L}\right)-H^{(1)}_0\left(\frac{\omega|x|}{c_T}\right)\right),\\
\end{equation}
where $\delta_{ij}$ is the Kronecker's delta, $\p_j$ denotes
$\p/\p x_j$, $ c_T=\sqrt{\mu}$, $c_L=\sqrt{\lambda+2\mu}$, and
$H^{(1)}_0$ is the Hankel function of the first kind and of order
$0$. See \cite{AH}  and \cite[Chap. 2]{KUPR}. For $\omega =0$, we
set $\mathbf{\Gamma}^0$ to be the Kelvin matrix of fundamental
solutions to the Lam\'e system, {\it i.e.},
 \begin{equation}\label{defgamma}
 \Gamma_{ij}^0(x)=
\ds  \frac{\gamma_1}{2\pi}\delta_{ij}\ln|x|-
  \frac{\gamma_2}{2\pi}\frac{x_ix_j}{|x|^2},
 \end{equation}
where
 \beq \label{gagamama}
 \gamma_1=\frac12\left(\frac1{\mu}+\frac1{2\mu+\lambda}\right) \quad {\rm and} \quad \gamma_2=\frac12
 \left(\frac1{\mu} - \frac1{2\mu+\lambda}\right).
 \eeq

The single and double layer potentials of the density function
$\mathbf{\varphi}\in L^2(\p \Om)^2$ associated with the Lam{\'e}
parameters $(\lambda, \mu)$ are defined by
\begin{align}
&\mathcal{S}^\omega\mathbf{\varphi}(x)=\ds\int_{\p
\Om}\mathbf{\Gamma}^\omega(x-y)\mathbf{\varphi}(y)d\sigma(y), \quad x \in \RR^2, \label{single}\\
&\mathcal{D}^\omega\mathbf{\varphi}(x)=\ds\int_{\p
\Om}\pd{}{\nu_y}\mathbf{\Gamma}^\omega(x-y)\mathbf{\varphi}(y)d\sigma(y),
\quad x \in \RR^2 \setminus \p\Om
\end{align}
where $\pd{}{\nu_y}$ denotes the conormal derivative with respect to the $y$-variables.

The following formulae give the jump relations
 obeyed by the double layer potential and by the conormal derivative of the single
layer potential:
\begin{align}
\pd{(\mathcal{S}^\omega\mathbf{\varphi})}{\nu}\Big|_{\pm}(x)&=\Big(
\pm
\frac{1}{2}I+(\mathcal{K}^\omega)^*\Big)\mathbf{\varphi}(x),~~\mbox{a.e.}~x\in\p
\Om, \label{2.7} \\
(\mathcal{D}^\omega\mathbf{\varphi})\Big|_{\pm}(x)&=\Big( \mp
\frac{1}{2}I+\mathcal{K}^\omega\Big)\mathbf{\varphi}(x),~~\mbox{a.e.}~x\in\p
\Om, \label{2.8}
\end{align}
where the subscripts $\pm$ indicate the limit from outside and inside $\Om$, respectively,
the singular integral operator $\mathcal{K}^\omega$ is defined by
\begin{align}
\mathcal{K}^\omega(x)=\mbox{p.v.}\int_{\p \Om}
\pd{\mathbf{\Gamma}^\omega(x-y)}{\nu_y}\mathbf{\varphi}(y)d\sigma(y),\label{kintegral}
\end{align}
and $(\mathcal{K}^\omega)^*$ is its $L^2$-adjoint, that is,
\begin{align}
(\mathcal{K}^\omega)^*(x)=\mbox{p.v.}\int_{\p \Om}
\pd{\mathbf{\Gamma}^\omega(x-y)}{\nu_x}\mathbf{\varphi}(y)d\sigma(y).
\end{align}
Here p.v. means the Cauchy principal value. See \cite{KUPR, DKV}.

It is proved in \cite{AAL,DKV} that $\mathcal{S}^\omega$ is
Fredholm as an operator from $L^2(\p \Om)^2$ into $H^1(\p \Om)^2$
and has index zero. It is also proved that the operators
$\pm\frac{1}{2}I+\mathcal{K}^\omega$ and
$\pm\frac{1}{2}I+(\mathcal{K}^\omega)^*$ are also Fredholm on
$L^2(\p \Om)^2$ with index $0$. We should emphasize that
$\mathcal{K}^\omega$ is not compact, even on bounded
$\mathcal{C}^\infty$-domains \cite{DKV}.

Let $\Psi$ be the vector space of all linear solutions to the
equation $\mathcal{L}^{\lambda,\mu}\mathbf{u}=0$ and
$\ds\pd{\mathbf{u}}{\nu}=0$ on $\p \Om$, or alternatively,
\begin{align*}
\Psi=\big\{\mathbf{\psi} : \p_i\psi_j+\p_j\psi_i=0, ~1\leq i,j
\leq 2 \big\}.
\end{align*}
Here $\psi_i$ for $i=1,2,$ denote the components of
$\mathbf{\psi}$. Define
\begin{align*}
L^2_\Psi(\p \Om)=\Big\{\mathbf{f}\in L^2(\p \Om)^2 : \ds\int_{\p
\Om} \mathbf{f}\cdot \mathbf{\psi} d\sigma=0 ~\mbox{for
all}~\mathbf{\psi} \in \Psi \Big\} ,
\end{align*}
a subspace of codimension $3$ in $L^2(\p \Om)^2$. In particular,
since $\Psi$ contains constant functions, we get
\begin{align*}
\ds\int_{\p \Om} \mathbf{f} d\sigma=0 \quad \mbox{for any}~
\mathbf{f}\in L^2_\Psi(\p \Om).
\end{align*}

We recall Green's formulae for the Lam{\'e} system, which can be
obtained by integration by parts. The first formula is
\begin{align}\label{greenfirst}
\int_{\p \Om}\mathbf{u}\cdot\pd{\mathbf{v}}{\nu}d\sigma=\int_{\Om}
\mathbf{u}\cdot\mathcal{L}^{\lambda,\mu}\mathbf{v}+\mathbf{E}(\mathbf{u},\mathbf{v}),
\end{align}
where $\mathbf{u}\in H^1(\Om)^2$, $\mathbf{v}\in
H^{\frac{3}{2}}(\Om)^2$, and
\begin{align} \label{E}
\mathbf{E}(\mathbf{u},\mathbf{v})= \int_\Om \lambda(\nabla\cdot\mathbf{u})(\nabla\cdot\mathbf{v})+\frac{\mu}{2}
(\nabla\mathbf{u}+\nabla\mathbf{u}^t)\cdot
(\nabla\mathbf{v}+\nabla\mathbf{v}^t).
\end{align}
Formula \eqnref{greenfirst} yields Green's second formula
\begin{align}\label{greensecond}
\int_{\p
\Om}\left(\mathbf{u}\cdot\pd{\mathbf{v}}{\nu}-\mathbf{v}\cdot\pd{\mathbf{u}}{\nu}\right)=\int_{\Om}
(\mathbf{u}\cdot\mathcal{L}^{\lambda,\mu}\mathbf{v}-\mathbf{v}\cdot\mathcal{L}^{\lambda,\mu}\mathbf{u}),
\end{align}
where $\mathbf{u},~ \mathbf{v}\in H^{\frac{3}{2}}(\Om)^2$.

Formula \eqnref{greensecond} shows that if $\mathbf{u}\in
H^{\frac{3}{2}}(\Om)^2$ satisfies
$\mathcal{L}^{\lambda,\mu}\mathbf{u}=0$ in $\Om$, then
$\ds\pd{\mathbf{u}}{\nu}\Big|_{\p \Om}\in L^2_\Psi(\p \Om)$.

\subsection{Radiation conditions and uniqueness}
Let us now formulate the {\it radiation condition} for the elastic
waves when $\mbox{Im }\omega\geq 0$ and $\omega \neq 0$.  See
\cite{AH, AAL, KUPR}. Any solution $\mathbf{u}$ to
$(\mathcal{L}^{\lambda,\mu}+\omega^2)\mathbf{u}=0$ admits the
decomposition,
 \beq \label{uupus}
 \mathbf{u}=\mathbf{u}^{(p)}+\mathbf{u}^{(s)},
 \eeq
where $\mathbf{u}^{(p)}$ and $\mathbf{u}^{(s)}$ are given by
\begin{align}
\mathbf{u}^{(p)} & =(k_T^2-k_L^2)^{-1}(\triangle + k_T^2)\mathbf{u},\\
\mathbf{u}^{(s)} & =(k_L^2-k_T^2)^{-1}(\triangle +
k_L^2)\mathbf{u},
\end{align}
with
 \beq \label{klkt}
 k_T=\frac{\omega}{c_T}=\frac{\omega}{\sqrt{\mu}}\quad \mbox{and}
 \quad k_L=\frac{\omega}{c_L}=\frac{\omega}{\sqrt{\lambda+2\mu}}.
 \eeq
Then $\mathbf{u}^{(p)}$ and $\mathbf{u}^{(s)}$ satisfy the
equations
 \beq \label{upuseqn}
 \begin{cases}
 (\triangle + k_T^2)\mathbf{u}^{(p)}=0, \quad & \rot\mathbf{u}^{(p)}=0,\\
 (\triangle + k_L^2)\mathbf{u}^{(s)}=0, \quad & \nabla\cdot\mathbf{u}^{(s)}=0.
 \end{cases}
 \eeq

We impose on $\mathbf{u}^{(p)}$ and $\mathbf{u}^{(s)}$ the
Sommerfeld radiation conditions for the solutions of the
corresponding Helmholtz equations by requiring that
 \beq \label{upusrad}
 \begin{cases}
 \ds \p_r\mathbf{u}^{(p)}(x)-\sqrt{-1}k_T\mathbf{u}^{(p)}(x)=o(r^{-\frac{1}{2}}),\\
 \nm
 \ds \p_r\mathbf{u}^{(s)}(x)-\sqrt{-1}k_L\mathbf{u}^{(s)}(x)=o(r^{-\frac{1}{2}}),
 \end{cases}
 \quad \mbox{as } r=|x|\rightarrow + \infty.
 \eeq
We say that $\mathbf{u}$ satisfies the radiation condition if it
allows the decomposition \eqnref{uupus} with $\mathbf{u}^{(p)}$
and $\mathbf{u}^{(s)}$ satisfying \eqnref{upuseqn} and
\eqnref{upusrad}. By a straightforward calculation one can see
that (non-periodic) single and double layer potentials,
$\Scal^\omega \phi$ and $\Dcal^\omega \phi$, satisfy the radiation
condition. See \cite{AAL,KUPR}.

We recall the following uniqueness results for  the exterior
problem \cite{KUPR}.
\begin{lem}\label{uniqueness1}
Let $\mathbf{u}$ be a solution to
$(\mathcal{L}^{\lambda,\mu}+\omega^2)\mathbf{u}=0$ in
$\mathbb{R}^2\setminus \overline{\Om}$ satisfying the radiation
condition. If either $\mathbf{u}=0$ or $\pd{\mathbf{u}}{\nu}=0$ on
$\p \Om$, then $\mathbf{u}$ is identically zero in
$\mathbb{R}^2\setminus \overline{\Om}$.
\end{lem}

Let $U$ be a bounded and connected open subset of $\mathbb{R}^2$
with the Lipschitz boundary such that $\overline{\Om} \subset U$.
Let $D_1 = U \setminus \overline{\Om}$ and
$D_2=(\mathbb{R}^2\setminus
\overline{\Om})\setminus\overline{D}_1$. Let $\mathcal{L}^{\tilde{\lambda},\tilde{\mu}}$ be the operator defined by \eqnref{lamedef}
and $\pd{}{\tilde{\nu}}$ be the corresponding conormal derivative. Consider the following
two-phase transmission problem:
\begin{equation}\label{twophase}
\begin{cases}
\mathcal{L}^{\tilde{\lambda},\tilde{\mu}}\mathbf{u}+\omega^2\mathbf{u}=0,
\quad & \mbox{in}~~ D_1,\\
\mathcal{L}^{\lambda,\mu}\mathbf{u}+\omega^2\mathbf{u}=0, \quad &
\mbox{in}~~ D_2,\\
\mathbf{u}\big |_+-\mathbf{u}\big |_-=0, \quad &\mbox{on}~~\p D_1\cap \p D_2,\\
\nm \ds\pd{\mathbf{u}}{\nu}\big
|_+-\pd{\mathbf{u}}{\tilde{\nu}}\big |_-=0, \quad &\mbox{on}~~\p
D_1\cap \p D_2,
\end{cases}
\end{equation}
with the radiation condition. The following uniqueness result
holds.
\begin{lem}\label{uniqueness2}
Let $\mathbf{u}$ be a solution to \eqnref{twophase} in
$\mathbb{R}^2\setminus \overline{\Om}$. If either $\mathbf{u}=0$
or $\pd{\mathbf{u}}{\nu}=0$ on $\p \Om$, $\mathbf{u}$ is
identically zero in $\mathbb{R}^2\setminus \overline{\Om}$.
\end{lem}

We note that the above two lemmas hold even when $\Om$ is an empty
set.

\subsection{Quasi-periodic Green's function}
In this section we collect some notation and well-known results
regarding quasi-periodic layer potentials for the Lam\'{e} system
in $\mathbb{R}^2$. We refer to \cite{mcphedran, poul, harris} for
the details.

We assume that the unit cell $Y=[0,1[^2$ is the periodic cell and
the quasi-momentum variable, denoted by $\alpha$, ranges over the
Brillouin zone $B= [0, 2 \pi[^2$. We introduce the two-dimensional
quasi-periodic Green's function $\mathbf{G}^{\alpha,\omega}$ for
$\omega \neq 0$, which satisfies
 \begin{equation} \label{fG}
 (\mathcal{L}^{\lambda,\mu} + \omega^2) \mathbf{G}^{\alpha,\omega} (x,y) =
 \sum_{n \in \ZZ^2} \delta(x-y - n) e^{\sqrt{-1} n\cdot \alpha}I.
 \end{equation}
Here we assume that $k_T, k_L \neq |2\pi n + \alpha|$ for all $n
\in \ZZ^2$ where $k_T$ and $k_L$ are given by \eqnref{klkt}.

A function $\mathbf{u}$ is said to be quasi-periodic or
$\alpha$-quasi-periodic if $e^{-\sqrt{-1}\alpha \cdot x}
\mathbf{u}$ is periodic. Using Poisson's summation formula, we
have
 $$
 \sum_{n \in \mathbb{Z}^2} \delta(x-y-n) e^{\sqrt{-1} n\cdot\alpha }I= \sum_{n \in \mathbb{Z}^2}
 e^{\sqrt{-1} (2 \pi n +\alpha) \cdot (x-y)}I.
 $$
We plug this equation into \eqnref{fG} and then take the Fourier
transform of both sides of \eqnref{fG} to obtain
 \begin{align*}
 \widehat{G}^{\alpha,\omega}_{ij}(\xi,y) = &(2\pi)^2\left\{ \frac{\delta_{ij}}{c_T^2}\frac{1}{k_T^2 -
\xi^2}+ \frac{\xi_i\xi_j}{\omega^2}\left(\frac{1}{k_L^2 -
\xi^2}-\frac{1}{k_T^2 - \xi^2}\right)\right\}\\
&\times\sum_{n\in \ZZ^2}e^{-\sqrt{-1}(2\pi n+ \alpha) \cdot
y}\delta(\xi+2\pi n+ \alpha) ,
 \end{align*}
where $\xi^2=\xi\cdot\xi$ and $~\widehat{}~$ denotes the Fourier
transform. Then taking the inverse Fourier transform, we can see
that the quasi-periodic Green's function
$\mathbf{G}^{\alpha,\omega}=(G^{\alpha,\omega}_{ij})$ can be
represented as a sum of augmented plane waves over the reciprocal
lattice:
\begin{align}
&G^{\alpha,\omega}_{ij}(x,y) =  \frac{\delta_{ij}}{c_T^2}\sum_{n\in
\ZZ^2}\frac{e^{\sqrt{-1}(2\pi n+ \alpha) \cdot (x-y)}}{k_T^2 -
|2\pi n
+\alpha|^2}\nonumber\\
&+ \frac{k_T^2-k_L^2}{\omega^2}\sum_{n\in
\ZZ^2}\frac{e^{\sqrt{-1}(2\pi n+ \alpha) \cdot (x-y)}(2\pi n+
\alpha)_i(2\pi n+ \alpha)_j}{(k_L^2 - |2\pi n + \alpha|^2)(k_T^2 -
|2\pi n + \alpha|^2)}.\label{reciplat}
\end{align}
Moreover, it can also be easily shown that
\begin{equation} \label{rephelm} \mathbf{G}^{\alpha,\omega} (x,y) =  \sum_{n \in \ZZ^2} \mathbf{\Gamma}^\omega(x- n -y) e^{\sqrt{-1}
n\cdot \alpha},
\end{equation} where $\mathbf{\Gamma}^\omega$ is the Green's matrix defined by \eqnref{kupradze}.

When $\omega=0$, we define $\mathbf{G}^{\alpha,0}$ by
 \beq \label{galpha0}
 G^{\alpha,0}_{ij}(x,y):= \frac{1}{\mu} \sum_{n\in \ZZ^2}e^{\sqrt{-1}
(2\pi n+ \alpha) \cdot (x-y)}\left( \frac{-\delta_{ij}}{ |2\pi n +
\alpha|^{2}}+ \frac{\lambda+\mu}{\lambda+2\mu} \frac{(2\pi n+
\alpha)_i(2\pi n+ \alpha)_j}{ |2\pi n + \alpha|^{4}}\right)
 \eeq
if $\alpha \neq 0$, while if $\alpha=0$, it is defined by
 \beq \label{galpha1}
 G^{0,0}_{ij}(x,y):= \frac{1}{\mu} \sum_{n \neq (0,0)}e^{\sqrt{-1}
 2\pi n \cdot (x-y)}\left( \frac{-\delta_{ij}}{ |2\pi n |^{2}}+ \frac{\lambda+\mu}{\lambda+2\mu}
 \frac{4\pi^2 n_i n_j}{ |2\pi n |^{4}}\right)
 \eeq

Then $\mathbf{G}^{\alpha,0}$ is quasi-periodic and satisfies
 \begin{align}
 \mathcal{L}^{\lambda,\mu}\mathbf{G}^{\alpha,0}(x,y) &=\sum_{n\in
 \ZZ^2}\delta(x-y-n)I  \quad \mbox{if } \alpha \neq 0, \label{g01} \\
 \mathcal{L}^{\lambda,\mu}\mathbf{G}^{0,0}(x,y) &=\sum_{n\in
 \ZZ^2}\delta(x-y-n)I-I. \label{g00}
 \end{align}
See \cite{indiana, book2} for the proof.

Let $D$ be a bounded domain in $\RR^2$ with a connected Lipschitz
boundary $\partial D$.  Let $\Scal^{\alpha,\omega}$ and
$\Dcal^{\alpha,\omega}$ be the quasi-periodic single and double
layer
 potentials
associated with $\mathbf{G}^{\alpha,\omega}$, that is, for a given
density $\vp \in L^2(\p D)^2$,
\begin{align*}
\Scal^{\alpha, \omega} \vp (x) & = \int_{\p D}
\mathbf{G}^{\alpha,\omega}(x,y) \vp (y) \, d\sigma (y), \quad x
\in \RR^2,
\\ \Dcal^{\alpha,\omega} \vp (x) & = \int_{\p D}
\pd{\mathbf{G}^{\alpha,\omega} (x,y)}{\nu_y} \vp (y) \, d\sigma
(y), \quad x \in \RR^2 \setminus \p D,
\end{align*}
where $\pd{}{\nu_y}$ denotes the conormal derivative with respect
to $y$. Then, $\Scal^{\alpha, \omega} \vp$ and
$\Dcal^{\alpha,\omega} \vp$ are solutions to  \begin{align*}
&(\mathcal{L}^{\lambda,\mu} + \omega^2) \mathbf{u}=0
\end{align*}
in $D$ and $Y\setminus\overline{D}$  and they are
$\alpha$-quasi-periodic.

The next  formulae give the jump relations obeyed by the double
layer potential and by the normal derivative of the single layer
potential on general Lipschitz domains:
\begin{align}
\pd{(\Scal^{\alpha,\omega} \vp)}{\nu} \bigg |_{\pm}(x) & = \bigg
(\pm \frac{1}{2} I + (\Kcal^{\alpha,\omega})^* \bigg ) \vp (x)
\quad \mbox{a.e. } x \in \p D \label{nuSp}, \\
 (\Dcal^{\alpha,\omega} \vp) \bigg|_{\pm} (x) & = \bigg(\mp \frac{1}{2}
I + \Kcal^{-\alpha,\omega} \bigg) \vp (x) \quad \mbox{a.e. } x \in
\p D, \label{doublejump-hp}
\end{align}
for $\vp \in L^2(\p D)^2$, where $\Kcal^{\alpha,\omega}$ is the
operator on $L^2(\p D)^2$ defined by
 \beq
 \Kcal^{\alpha,\omega} \vp
 (x) = \mbox{p.v.} \int_{\p D}
 \pd{\mathbf{G}^{-\alpha,\omega}(x,y)}{\nu_y} \vp (y) d \sigma
 (y),
 \eeq
 and $(\Kcal^{\alpha,\omega})^*$ is given by
 \beq \label{kcalao}
 (\Kcal^{\alpha,\omega})^* \vp (x) = \mbox{p.v.} \int_{\p D}
 \pd{\mathbf{G}^{\alpha,\omega}(x,y)}{\nu_x} \vp (y) d \sigma (y).
 \eeq
Note that $(\Kcal^{\alpha,\omega})^*$ is the $L^2$-adjoint of
$\Kcal^{\alpha,\omega}$. The formulae \eqnref{nuSp} and \eqnref{doublejump-hp} hold because
$\mathbf{G}^{\alpha,\omega}(x,y)$ has the same kind of singulary at $x=y$ as that of $\mathbf{\Gamma}(x-y)$.

The following lemma will be of use in later sections
\begin{lem}
For any constant vector $\phi$
 \beq\label{d001}
 (\frac{1}{2} I + \Kcal^{0,0})\phi = |Y\setminus D| \phi \quad
 \mbox{on } \p D,
 \eeq
and for any $\psi \in L^2(\p D)^2$
 \beq \label{d002}
 \int_{\p D} \left(\frac{1}{2}I+ (\Kcal^{0,0})^*\right)\psi = |Y\setminus
 D| \int_{\p D} \psi.
 \eeq
\end{lem}
\pf By Green's theorem and \eqnref{g00} we have
 $$
 \Dcal^{0,0}\phi (x)= \int_{D}
 \mathcal{L}^{\lambda,\mu}\mathbf{G}^{0,0}(x,y)\phi dy= \phi -
 \int_D \phi,
 $$
and hence \eqnref{d001} follows since $(\frac{1}{2} I +
\Kcal^{0,0})\phi = \Dcal^{0,0}\phi|_{-}$.

The identity \eqnref{d002} immediately follows from \eqnref{d001}. In fact, for any constant vector $\phi$, we have
 \begin{align*}
 \int_{\p D} \phi \cdot \left(\frac{1}{2}I+ (\Kcal^{0,0})^*\right)\psi &=
 \int_{\p D} \left(\frac{1}{2}I+ \Kcal^{0,0} \right) \phi \cdot \psi \\
 & = |Y\setminus D| \int_{\p D} \phi \cdot \psi.
 \end{align*}
Thus \eqnref{d002} follows.
\qed

\section{Asymptotic behavior of phononic bands} The phononic
crystal we consider in this paper is a homogeneous elastic medium
of Lam\'{e} constants $\lambda$, $\mu$ which contains an array of
arbitrary-shaped inclusions $\Om=\cup_{n\in \mathbb{Z}^2}(D+n)$
which is periodic with respect to the lattice $\mathbb{Z}^2$.
These inclusions have Lam\'{e} constants $\tilde{\lambda}$,
$\tilde{\mu}$. Let $Y =]0,1[^2$ denote the fundamental period
cell. For each quasi-momentum $\alpha\in[0,2\pi[^2$, let
$\sigma_\alpha(D)$ be the (discrete) spectrum of the operator
defined by \eqnref{periodicop} with the condition that  $e^{-
\sqrt{-1} \alpha \cdot x}\mathbf{u}$ is periodic. In other words,
$\sigma_\alpha (D)$ is the spectrum of the problem
\begin{equation}\label{transprime3}
\begin{cases}
\mathcal{L}^{\lambda,\mu}\mathbf{u}+\omega^2\mathbf{u}=0, \quad &
\mbox{in}~~ Y\setminus
\overline{D},\\
\mathcal{L}^{\tilde{\lambda},\tilde{\mu}}\mathbf{u}+\omega^2\mathbf{u}=0,
\quad & \mbox{in}~~ D,\\
\nm
\mathbf{u}\big |_+-\mathbf{u}\big |_-=0, \quad &\mbox{on}~~\p D,\\
\nm \ds\pd{\mathbf{u}}{\nu}\big
|_+-\pd{\mathbf{u}}{\tilde{\nu}}\big |_-=0, \quad &\mbox{on}~~\p
D,\\
\nm e^{- \sqrt{-1} \alpha \cdot x}\mathbf{u} \mbox{ is periodic}.
\end{cases}\end{equation}
Here $\mathcal{L}^{\tlambda,\tmu}$ is the  elastostatic system
corresponding to the Lam{\'e} constants $\tlambda$ and $\tmu$ and
${\partial}/{\partial \tnu}$ denotes the corresponding conormal
derivative.

By the standard Floquet theory, the spectrum of
\eqnref{transprime3} has the band structure given by
 \beq
 \ds \bigcup_{\alpha \in [0,2\pi[^2} \sigma_\alpha(D).
 \eeq
The main objective of this section is to investigate the behavior
of $\sigma_\alpha(D)$ as $\mu \to +\infty$.

\subsection{Integral representation of quasi-periodic solutions}

In this section, we obtain the integral representation formula for
the solution to \eqnref{transprime3}. We denote by
$\tilde{\Scal}^{\omega}, \tilde{\Dcal}^{\omega}$, and
$\tilde{\Kcal}^{\omega}$ the layer potentials associated with the
Lam{\'e} parameters $(\tlambda, \tmu)$.

We first prove the following lemma.
\begin{lem} \label{ortho}
Suppose that $\omega^2$ is not an eigenvalue for
$-\mathcal{L}^{\lambda,\mu}$ in $D$ with the Dirichlet boundary
condition on $\partial D$. Let $\mathbf{u}$ be a solution to
\eqnref{transprime3}. Then  we have
$$\mathbf{u}|_{\partial
D}\perp\ker\tilde{\Scal}^{\omega} \quad\mbox{and} \quad
\mathbf{u}|_{\partial D}\perp\ker(\Scal^{\alpha,\omega})^*.$$ Here
$\tilde{\Scal}^{\omega}$ and $\Scal^{\alpha,\omega}$ are
considered to be operators on $L^2(\p D)^2$.
\end{lem}

\pf We first observe that, since
$(\mathcal{L}^{\tilde{\lambda},\tilde{\mu}}+\omega^2)\mathbf{u}=0$
in $D$, we have
 \beq \label{greenid}
 \mathbf{u}(x)=\tilde{\Dcal}^{\omega}\left(\mathbf{u}|_{\partial
 D}\right)(x)-\tilde{\Scal}^{\omega}\left(\frac{\partial
 \mathbf{u}}{\partial\tilde{\nu}}\Big|_-\right)(x), \quad x  \in D,
 \eeq
and consequently by \eqnref{2.8}
 \beq \label{3.1.1}
\ds\frac{1}{2}\mathbf{u}|_{\partial
D}=\tilde{\Kcal}^{\omega}\left(\mathbf{u}|_{\partial
D}\right)-\tilde{\Scal}^{\omega}\left(\frac{\partial
\mathbf{u}}{\partial\tilde{\nu}}\Big|_-\right).
 \eeq

Let $\phi\in\mbox{Ker}(\tilde{\Scal}^{\omega})$, {\it i.e.},
$\tilde{\Scal}^{\omega}\phi=0$ on $\p D$. By Lemma
\ref{uniqueness1}, we have $\tilde{\Scal}^{\omega}\phi=0$ in
$\RR^2\setminus D$ and hence
$\frac{1}{2}\phi+(\tilde{\Kcal}^{\omega}\big)^*\phi=0$ by
\eqnref{2.7}. Then we have from \eqnref{3.1.1}
\begin{align*}
\frac{1}{2} \langle \mathbf{u}|_{\partial D},\phi \rangle & =\left
\langle \tilde{\Kcal}^{\omega}\left(\mathbf{u}|_{\partial
D}\right),\phi\right\rangle
-\left\langle \tilde{\Scal}^{\omega}\left(\frac{\partial \mathbf{u}}{\partial\tilde{\nu}}\Big|_-\right),\phi\right\rangle\\
&=\left\langle \mathbf{u}|_{\partial
D},(\tilde{\Kcal}^{\omega})^*\phi\right\rangle
-\left\langle \frac{\partial \mathbf{u}}{\partial\tilde{\nu}}\Big|_-,\tilde{\Scal}^{\omega}\phi\right\rangle \\
&= -\frac{1}{2}\left\langle \mathbf{u}|_{\partial
D},\phi\right\rangle,
\end{align*}
which implies $\langle\mathbf{u}|_{\partial D},\phi\rangle=0$, and
hence $\mathbf{u}|_{\partial D}\perp\ker\tilde{\Scal}^{\omega}$.

Observe that if $\mathbf{u}$ is $\alpha$-quasi-periodic, then
\begin{align*}
\Dcal_Y^{\alpha,\omega}\left(\mathbf{u}|_{\partial Y}\right)=0
\quad \mbox{and} \quad \Scal_Y^{\alpha,\omega}\left(\frac{\partial
\mathbf{u}}{\partial\nu}\Big|_+\right)=0 \quad \mbox{on } \p Y,
\end{align*}
where $\Dcal_Y^{\alpha,\omega}$ and $\Scal_Y^{\alpha,\omega}$ are
the ($\alpha$-quasi-periodic) double and single layer potentials
on $\p Y$. Thus we have
$$
\mathbf{u}(x)=-\Dcal^{\alpha,\omega}\left(\mathbf{u}|_{\partial
D}\right)(x)+\Scal^{\alpha,\omega}\left(\frac{\partial
\mathbf{u}}{\partial\nu}\Big|_+\right)(x), \quad x  \in Y\setminus
\overline{D},
$$
and consequently,
$$
\frac{1}{2}\mathbf{u}|_{\partial
D}=-\Kcal^{-\alpha,\omega}\left(\mathbf{u}|_{\partial
D}\right)+\Scal^{\alpha,\omega}\left(\frac{\partial
\mathbf{u}}{\partial\nu}\Big|_+\right).
$$

Let $\phi\in\ker(\Scal^{\alpha,\omega})^*$. Since
$(\Scal^{\alpha,\omega})^*=\Scal^{-\alpha,\omega}$, we have
$$\Scal^{-\alpha,\omega} \phi=0\quad \mbox{on}~ \p D.$$
Since $\omega^2$ is not a Dirichlet eigenvalue of
$-\mathcal{L}^{\lambda,\mu}$ in $D$, we immediately deduce that
$$\Scal^{-\alpha,\omega} \phi=0\quad \mbox{in}~  D,$$ and hence
$$-\frac{1}{2}\phi+(\Kcal^{-\alpha,\omega}\big)^*\phi=0 \quad \mbox{on} ~\p D.$$
Therefore, we get
 \begin{align*}
 \frac{1}{2} \langle \mathbf{u}|_{\partial D},\phi \rangle &
=-\left\langle \Kcal^{-\alpha,\omega}\left(\mathbf{u}|_{\partial
D}\right),\phi\right\rangle
+\left\langle \Scal^{\alpha,\omega}\left(\frac{\partial \mathbf{u}}{\partial\nu}\Big|_+\right),\phi\right\rangle \\
&=-\left\langle \mathbf{u}|_{\partial
D},(\Kcal^{-\alpha,\omega})^*\phi\right\rangle
+\left\langle \frac{\partial \mathbf{u}}{\partial\nu}\Big|_+,\Scal^{-\alpha,\omega}\phi\right\rangle \\
&= -\frac{1}{2}\left\langle \mathbf{u}|_{\partial
D},\phi\right\rangle ,
\end{align*}
which implies $\left\langle \mathbf{u}|_{\partial
D},\phi\right\rangle =0$. This completes the proof. \qed

We now establish a representation formula for solutions of
\eqnref{transprime3}.

\begin{thm} \label{repthm}
Suppose that $\omega^2$ is not an eigenvalue for
$-\mathcal{L}^{\lambda,\mu}$ in $D$ with the Dirichlet boundary
condition on $\partial D$. Then, for any solution $\mathbf{u}$ of
\eqnref{transprime3}, there exists one and only one pair
$(\phi,\psi)\in L^2(\partial D)^2\times L^2(\partial D)^2$ such
that
 \begin{equation}\label{repfor}
 \mathbf{u}(x)= \left\{ \begin{array}{ll} \tilde{\Scal}^{\omega} \phi (x),
 \quad & x \in D,\\
 \nm
 \Scal^{\alpha, \omega} \psi (x),
 \quad & x \in Y \setminus \overline{D}. \end{array} \right.
\end{equation}
Moreover, $(\phi, \psi)$ satisfies
\begin{equation}\label{inteq2}
\begin{cases}
\ds \tilde{\Scal}^{\omega} \phi - \Scal^{\alpha,
\omega} \psi = 0 \quad & \mbox{on } \partial D,\\
\nm \ds \Big(\frac{1}{2}I- (\tilde{\Kcal}^{\omega})^* \Big )\phi +
\Big(\frac{1}{2}I+ (\Kcal^{\alpha, \omega})^* \Big )\psi =0 \quad
& \mbox{on }
\partial D,
\end{cases}
\end{equation}
and the mapping $\mathbf{u}\longmapsto (\phi,\psi)$ from solutions
of \eqnref{transprime3} in $H^1(Y)^2$ to solutions to the system
of integral equations \eqnref{inteq2} in $L^2(\partial D)^2\times
L^2(\partial D)^2$ is a one-to-one correspondence.
\end{thm}
\pf We first note that the problem of finding $(\phi,\psi)$
satisfying  \eqnref{repfor} and \eqnref{inteq2} is equivalent to
solving the following two systems of equations:
\begin{equation}
  \label{eq:phi}
\begin{cases}
\tilde{\Scal}^{\omega}\phi=\mathbf{u}|_{\partial D}& \text{ on }\partial D,\\
\nm
\ds\big(-\frac{1}{2}I+(\tilde{\Kcal}^{\omega})^*\big)\phi=\frac{\partial
\mathbf{u}}{\partial\tilde{\nu}}\Big|_{-}& \text{ on }\partial D,
\end{cases}
\end{equation}
and
\begin{equation}
  \label{eq:psi}
\begin{cases}
\Scal^{\alpha,\omega}\psi=\mathbf{u}|_{\partial D}& \text{ on }\partial D,\\
\nm
\ds\big(\frac{1}{2}I+(\Kcal^{\alpha,\omega})^*\big)\psi=\frac{\partial
\mathbf{u}}{\partial\nu}\Big|_{+}& \text{ on }\partial D.
\end{cases}
\end{equation}

In order to find $\phi$ satisfying \eqnref{eq:phi}, it suffices to
find $\phi$ satisfying $\tilde{\Scal}^{\omega}\phi=\mathbf{u}$ in
$D$. Since $\tilde{\Scal}^{\omega}$ is self-adjoint Fredholm
operator on $L^2(\p D)^2$ with index $0$ \cite{AAL}, it follows
from Lemma \ref{ortho} that there exists $\phi_0\in L^2(\partial
D)^2$ such that
 \beq \label{sphiu}
\tilde{\Scal}^{\omega}\phi_0=\mathbf{u}|_{\partial D} \quad\text{
on }\partial D.
 \eeq
Observe that if $\omega \neq 0$, then the solution to the
Dirichlet problem for $\mathcal{L}^{\lambda,\mu}+\omega^2$ may not
be unique, and hence \eqnref{sphiu} does not imply
$\tilde{\Scal}^{\omega}\phi_0=\mathbf{u}$ in $D$. However, since
$(\mathcal{L}^{\lambda,\mu}+\omega^2) (\mathbf{u}-
\tilde{\Scal}^{\omega}\phi_0)=0$ in $D$, we get by Green's formula
\begin{align*}
\mathbf{u}-\tilde{\Scal}^{\omega}\phi_0=-\tilde{\Scal}^{\omega}
\left[\pd{}{\tilde{\nu}}\left(
\mathbf{u}-\tilde{\Scal}^{\omega}\phi_0 \right)\Big|_- \right]
\quad\text{ in }D,
\end{align*}
and therefore,
\begin{align*}
\mathbf{u}=\tilde{\Scal}^{\omega}\left[\phi_0-\pd{}{\tilde{\nu}}\left(
\mathbf{u}-\tilde{\Scal}^{\omega}\phi_0
\right)\Big|_-\right]\quad\text{ in }D.
\end{align*}

To prove the uniqueness of $\phi$ satisfying \eqnref{eq:phi}, it
suffices to show that the solution to
\begin{equation*}
\begin{cases}
\tilde{\Scal}^{\omega}\phi=0& \text{ on }\partial D,\\
\nm \ds\big(-\frac{1}{2}I+(\tilde{\Kcal}^{\omega})^*\big)\phi=0&
\text{ on }\partial D,
\end{cases}
\end{equation*} is zero.  By the first equation in \eqnref{eq:phi} and Lemma
\ref{uniqueness1}, $\tilde{\Scal}^{\omega}\phi=0$ in
$\RR^2\setminus\overline{D}$ and hence
$\phi=\ds\pd{}{\tilde{\nu}}\tilde{\Scal}^{\omega}\phi\Big|_+-\pd{}{\tilde{\nu}}\tilde{\Scal}^{\omega}\phi\Big|_-=0$.

Similarly, we can show existence and uniqueness of $\psi$
satisfying
\begin{align*}
\mathbf{u}=\Scal^{\alpha, \omega} \psi \quad\mbox{in}~
Y\setminus\overline{D},
\end{align*}
which yields \eqnref{eq:psi}. This completes the proof. \qed

Let $\Acal^{\alpha,\omega}$ be the operator-valued function of
$\omega$ defined by
 \begin{equation} \label{aalpha}
 \Acal^{\alpha,\omega} :=
 \left(\begin{array}{cc}
\ds \tilde{\Scal}^{\omega} & - \Scal^{\alpha,
\omega} \\
\nm \ds \frac{1}{2}I- (\tilde{\Kcal}^{\omega})^* &
 \ds\frac{1}{2}I+ (\Kcal^{\alpha, \omega})^*
  \end{array} \right).
 \end{equation}
By Theorem \ref{repthm}, $\omega^2$ is an eigenvalue corresponding
to quasi-momentum $\alpha$ if and only if $\omega$ is a
characteristic value of $\Acal^{\alpha,\omega}$. Consequently, we
have now a new way of computing the spectrum of
\eqnref{transprime3} by examining the characteristic values of
$\Acal^{\alpha,\omega}$. Based on Muller's method \cite{muller}
for finding complex roots of scalar equations, a boundary element
method similar to the one developed in \cite{cheng1,cheng2} can be
designed for computing phononic band gaps.

\subsection{Full asymptotic expansions}

Expanding the operator-valued functions $\Acal^{\alpha,\omega}$ in
terms of $\mu$ as $\mu\to + \infty$, we can calculate asymptotic
expressions of their characteristic values with the help of the
generalized Rouch{\'e}'s theorem, and it is what we do in this
subsection.

We begin with the following asymptotic expansion of
$G^{\alpha,\omega}_{ij}(x,y)$ in \eqnref{reciplat}.

\begin{lem}
Let $\tau_l=1-\left(\frac{c_T}{c_L}\right)^{2l}$. As $\mu \to + \infty$,
\begin{align}\label{gexp1}
G^{\alpha,\omega}_{ij} (x,y) &= \sum_{l=1}^{+
\infty}\frac{\omega^{2(l-1)}}{\mu^l}\sum_{n\in \ZZ^2}e^{\sqrt{-1}
(2\pi n+ \alpha) \cdot (x-y)}\left( \frac{-\delta_{ij}}{ |2\pi n +
\alpha|^{2l}}+\tau_l\frac{(2\pi n+ \alpha)_i(2\pi n+ \alpha)_j}{
|2\pi n + \alpha|^{2(l+1)}}\right),
\end{align}
for fixed $\alpha \neq 0$, while for $\alpha=0$,
\begin{align}\label{gexp2}
G^{0,\omega}_{ij} (x,y) =
\frac{\delta_{ij}}{\omega^2}+\sum_{l=1}^{+
\infty}\frac{\omega^{2(l-1)}}{(2\pi)^{2l}\mu^l}\sum_{n\in
\ZZ^2\setminus \{0\}}e^{2\pi \sqrt{-1} n \cdot (x-y)}\left(
-\frac{\delta_{ij}}{ |n|^{2l}}+\tau_l\frac{n_in_j}{
|n|^{2(l+1)}}\right).
\end{align}
\end{lem}

Derivation of \eqnref{gexp1} and \eqnref{gexp2} are
straightforward. In fact, since
 $$
 \frac{1}{k_T^2 - |2\pi n +\alpha|^2} = \frac{1}{\frac{\omega^2}{\mu} - |2\pi n +\alpha|^2}
 = -\sum_{k=0}^\infty \frac{\omega^{2k}}{\mu^k |2\pi n
 +\alpha|^{2(k+1)}},
 $$
one immediately obtains \eqnref{gexp1} and \eqnref{gexp2}.

We can write \eqnref{gexp1} and \eqnref{gexp2} as
 \beq
 \mathbf{G}^{\alpha,\omega}
 (x,y)=\sum_{l=1}^{+\infty}\frac{1}{\mu^l}\mathbf{G}_l^{\alpha,\omega}(x,y),
 \eeq
and
 \beq
 \mathbf{G}^{0,\omega} (x,y) =
 \frac{1}{\omega^2} {I}+\sum_{l=1}^{+\infty}\frac{1}{\mu^l}\mathbf{G}^{0,\omega}_l(x,y),
 \eeq
where the definitions of $\mathbf{G}_l^{\alpha,\omega}(x,y)$ and
$\mathbf{G}^{0,\omega}_l(x,y)$ are obvious from \eqnref{gexp1} and
\eqnref{gexp2}. We note that $\mathbf{G}_l^{\alpha,\omega}(x,y)$
and $\mathbf{G}^{0,\omega}_l(x,y)$ are dependent upon $\mu$
because of the factor $\tau_l$. However, since $|\tau_l| \le C$
for some constant $C$ independent of $\mu$ and $l$, it will not
affect analysis to follow. We also note that
$\mathbf{G}^{\alpha,\omega}_1(x,y)$ is independent of $\omega$ and
 \beq \label{314}
 \mathbf{G}^{\alpha,\omega}_1(x,y)= \mu
 \mathbf{G}^{\alpha,0}(x,y),
 \eeq
where $\mathbf{G}^{\alpha,0}(x,y)$ is the quasi-periodic Green
function defined in \eqnref{galpha0}.

Denote by $\Scal_{l}^{\alpha,\omega}$ and
$(\Kcal_{l}^{\alpha,\omega})^*$, for $l \geq 1$ and $\alpha \in
[0,2\pi[^2,$ the single layer potential and the boundary integral
operator associated with the kernel
$\mathbf{G}_l^{\alpha,\omega}(x,y)$ as defined in \eqnref{kcalao}
so that
 \begin{equation} \label{sumform}
 \Scal^{\alpha,\omega} = \sum_{l=1}^{+\infty}\frac{1}{\mu^l} \Scal_{l}^{\alpha,\omega} \quad
 \mbox{and} \quad (\Kcal^{\alpha,\omega})^* = \sum_{l=1}^{+\infty}\frac{1}{\mu^l} (\Kcal_{l}^{\alpha,\omega})^*.
 \end{equation}

\begin{lem} \label{invert}
The operator $\frac{1}{2} I + (\Kcal^{\alpha,0})^* : L^2(\partial
D)^2 \to  L^2(\partial D)^2$ is invertible.
\end{lem}

Before proving Lemma \ref{invert}, let us make a note of the
following simple fact: If $\mathbf{u}$ and $\mathbf{v}$ are
$\alpha$-quasi-periodic, then
 \begin{equation} \label{intpart}
 \int_{\p Y} \pd{\mathbf{u}}{\nu}\cdot\overline{\mathbf{v}} d\sigma =0.
 \end{equation}
To prove this, it is enough to see that
 \begin{align*}
 &\int_{\p Y}  \pd{\mathbf{u}}{\nu}\cdot\overline{\mathbf{v}} = \int_{\p Y}  \pd{(
 e^{-i\alpha\cdot x}\mathbf{u})}{\nu} \cdot \overline{e^{-i\alpha\cdot x}\mathbf{v}}\\
 &+i\int_{\p Y} \Big[  \lambda\alpha\cdot(e^{-i\alpha\cdot x}\mathbf{u})N+\mu\begin{pmatrix}2\alpha_1N_1+\alpha_2N_2 & \alpha_1N_2\\
 \alpha_2N_1 & \alpha_1N_1+2\alpha_2N_2 \end{pmatrix}(e^{-i\alpha\cdot x}\mathbf{u})   \Big]
 \cdot \overline{e^{-i\alpha\cdot x}\mathbf{v}}.
\end{align*}
Here $i= \sqrt{-1}$ and $N$ is the outward unit normal to $Y$.
Then the integrands over the opposite sides of $\p Y$ have the
same absolute values with different signs and therefore the integration over $\p Y$ is zero.\\

\noindent{\sl Proof of Lemma \ref{invert}}. For $\alpha\ne 0$, we
show injectivity of $\frac{1}{2} I + (\Kcal^{\alpha,0})^*$. Then
from the Fredholm alternative, the result follows. Suppose $\phi
\in L^2(\partial D)^2$ satisfies \begin{align*} (\frac{1}{2} I +
(\Kcal^{\alpha, 0})^*) \phi =0\quad\mbox{on}~\partial
D.\end{align*} Then by \eqnref{nuSp}, $\mathbf{u}:=
\Scal^{\alpha,0} \phi$ satisfies
\begin{align*}
\begin{cases}
\mathcal{L}^{\lambda,\mu} \mathbf{u}=0 \quad \mbox{in}~ Y
\setminus
\overline{D},\\
\ds\pd{\mathbf{u}}{\nu}\Big|_+ =0\quad \mbox{on}~\p D,\\
\mathbf{u}~\mbox{is}~ \alpha\mbox{-quasi-periodic}.
\end{cases}
\end{align*}
Therefore, it follows from \eqnref{intpart} that
 $$
 \int_{Y\setminus D}\left( \lambda|\nabla\cdot \mathbf{u}|^2+\frac{\mu}{2}|\nabla \mathbf{u}+\nabla \mathbf{u}^t|^2\right) = \int_{\p Y} \pd{\mathbf{u}}{\nu}\cdot\overline{\mathbf{u}}
 - \int_{\p D} \pd{\mathbf{u}}{\nu} \bigg|_+\cdot\overline{\mathbf{u}}=0.
 $$
Thus, $\mathbf{u}$ is constant in $Y\setminus \overline D$, and
hence in $D$. Hence, we get $$\phi =\ds
\pd{\mathbf{u}}{\nu}\Big|_+ - \pd{\mathbf{u}}{\nu}\Big|_-=0.$$

For the periodic case $\alpha=0$, we show the injectivity of
$\frac{1}{2} I + \Kcal^{0,0}$. Let $\phi \in L^2(\partial D)^2$
satisfying $(\frac{1}{2} I + \Kcal^{0, 0}) \phi =0$ on $\partial
D$. Then $\mathbf{u}:= \Dcal^{0,0} \phi$ satisfies
\begin{align*}
\begin{cases}
\mathcal{L}^{\lambda,\mu} \mathbf{u}=0 \quad & \mbox{in}~ D,\\
\mathbf{u}|_- =0 \quad & \mbox{on}~ \p D,  \\
\end{cases}
\end{align*}
and therefore $\mathbf{u}=0$ in $D$. Furthermore, if $(\frac{1}{2}
I + \Kcal^{0, 0}) \phi =0$, we can show that  $\phi\in H^1(\p
D)^2$ and $\ds\pd{(\Dcal^{0,0}
\phi)}{\nu}\Big|_+=\ds\pd{(\Dcal^{0,0} \phi)}{\nu}\Big|_-$. See
\cite{AAL} for  the details. Then we have
\begin{align*}
\begin{cases}
\mathcal{L}^{\lambda,\mu} \mathbf{u}=0 \quad \mbox{in}~ Y
\setminus
\overline{D},\\
\ds\pd{\mathbf{u}}{\nu}\Big|_+ =0\quad \mbox{on}~\p D,\\
\mathbf{u}~\mbox{is}~ \mbox{periodic}.
\end{cases}
\end{align*}
Therefore, it follows that
 $$
 \int_{Y\setminus D}\left( \lambda|\nabla\cdot \mathbf{u}|^2+\frac{\mu}{2}|\nabla \mathbf{u}+\nabla \mathbf{u}^t|^2\right) = \int_{\p Y} \pd{\mathbf{u}}{\nu}\cdot\overline{\mathbf{u}}
 - \int_{\p D} \pd{\mathbf{u}}{\nu} \bigg|_+\cdot \overline{\mathbf{u}}=0.
 $$
Thus, $\mathbf{u}$ is constant in $Y\setminus \overline D$, and
hence  $\phi =\mathbf{u}|_- - \mathbf{u}|_+$ is constant. By
\eqnref{d001}, we obtain that
\begin{align*}
0= (\frac{1}{2} I + \Kcal^{0,0})\phi =|Y\setminus D|\phi,
\end{align*}
which implies that $\phi$ must be zero. This completes the proof.
\qed

We now derive complete asymptotic expansions of eigenvalues as
$\mu\to + \infty$. We deal with three cases separately: $\alpha
\neq 0$ (not of order $O(\frac{1}{\sqrt{\mu}})$), $\alpha =0$, and
$|\alpha|$ is of order $O(\frac{1}{\sqrt{\mu}})$

\subsubsection{{The case $\alpha \neq 0$.}}

The following lemma, which is an immediate consequence of
\eqnref{sumform}, gives a complete asymptotic expansion of
$\Acal^{\alpha,\omega}$ defined in \eqnref{aalpha} as $\mu\to +
\infty$.

\begin{lem} \label{asympta} Suppose $\alpha \neq 0$.
Let
 \beq \label{317}
 \Acal^{\alpha,\omega}_0= \left(\begin{array}{cc}
 \tilde{\Scal}^{\omega}  & 0
 \\ \nm \ds \frac{1}{2}I- (\tilde{\Kcal}^{\omega})^* & \ds
 \frac{1}{2}I+ (\Kcal^{\alpha,0})^*
 \end{array} \right),
 \eeq
and, for $l\geq 1$,
 \beq \label{318}
 \Acal^{\alpha,\omega}_l=
 \left(\begin{array}{cc} 0 & - \Scal_{l}^{\alpha,\omega} \\
 \nm  0 & \frac{1}{\mu}(\Kcal_{l+1}^{\alpha,\omega})^*
 \end{array} \right).
 \eeq
Then we have
\begin{equation}\label{03p}
\Acal^{\alpha,\omega} = \Acal^{\alpha,\omega}_0
+\sum_{l=1}^{+\infty} \frac{1}{\mu^l}\Acal^{\alpha,\omega}_l.
\end{equation}
All the operators are defined on $L^2(\p D)^2\times L^2(\p D)^2$.
\end{lem}

Note that it is just for convenience that there is $1/\mu$ in the
definition of $\Acal^{\alpha,\omega}_l$. This of course does not
affect any of our asymptotic results.
\begin{lem} \label{chavalue}
Suppose $\alpha \neq 0$. Then the followings are equivalent:
\begin{itemize}
\item[\rm{(i)}] $\omega_0^\alpha \in \RR$ is a
characteristic value of $\Acal^{\alpha,\omega}_0$,
\item[\rm{(ii)}] $\omega_0^\alpha \in \RR$ is a
characteristic value of $\widetilde{\Scal}^{\omega}$,
\item[\rm{(iii)}] $(\omega_0^\alpha)^2$ is  an eigenvalue of
$-\mathcal{L}^{\tilde{\lambda},\tilde{\mu}}$ in $D$ with the
Dirichlet boundary condition.
\end{itemize}
Moreover if $\mathbf{u}$ is an eigenfunction of $-\mathcal{L}^{\tilde{\lambda},\tilde{\mu}}$ in $D$ with the
Dirichlet boundary condition, then $\vp := \pd{\mathbf{u}}{\nu}|_{-}$ is an characteristic function of $\widetilde{\Scal}^{\omega}$.
Conversely, if $\vp$ is an characteristic function of $\widetilde{\Scal}^{\omega}$, then $\mathbf{u}:=-\widetilde{\Scal}^{\omega} (\vp)$
is an eigenfunction of $-\mathcal{L}^{\tilde{\lambda},\tilde{\mu}}$ in $D$ with the
Dirichlet boundary condition.
\end{lem}
\pf By Lemma \ref{invert}, $\frac{1}{2}I+ (\Kcal^{\alpha,0})^*$ is
invertible. Thus characteristic values of
$\Acal^{\alpha,\omega}_0$ coincide with those of
$\widetilde{\Scal}^{\omega}$. On the other hand, the Green identity \eqnref{greenid} shows that
the characteristic values of $\widetilde{\Scal}^{\omega}$  are exactly eigenvalues of
$-\mathcal{L}^{\tilde{\lambda},\tilde{\mu}}$ in $D$ with the
Dirichlet boundary condition. The last statements of Lemma \ref{chavalue} also follow from \eqnref{greenid}.
\qed

\begin{lem}\label{rankone}
Every eigenvector of $\tilde{\Scal}^{\omega}$ has rank one.
\end{lem}
\pf Let $\phi$ be an eigenvector of $\tilde{\Scal}^{\omega}$
associated with the characteristic value $\omega_0$, {\it i.e.},
$\tilde{\Scal}^{\omega_0}\phi=0$ on $\p D$. Suppose that there
exists $\phi^\omega$, holomorphic in a neighborhood of $\omega_0$
as a function of $\omega$, such that $\phi^{\omega_0}=\phi$ and
\begin{align*}
\tilde{\Scal}^{\omega}\phi^\omega=(\omega^2-\omega_0^2)\psi^\omega
\end{align*}
for some $\psi^\omega$. Let
$\mathbf{u}^\omega(x):=\tilde{\Scal}^{\omega}\phi^\omega(x)$,
$x\in D$. Then $\mathbf{u}^\omega$ satisfies
\begin{align*}
\begin{cases}
(\mathcal{L}^{\tilde{\lambda},\tilde{\mu}}+\omega^2)\mathbf{u}^\omega=0\quad \mbox{in}~ D,\\
\mathbf{u}^\omega=(\omega^2-\omega_0^2)\psi^\omega\quad\mbox{on}~\p
D.
\end{cases}
\end{align*}
By Green's formula, we have
\begin{align*}
(\omega^2-\omega_0^2)\int_D \mathbf{u}^\omega
\cdot\overline{\mathbf{u}^{\omega_0}}&=\int_D \mathbf{u}^\omega
\cdot\overline{\mathcal{L}^{\tilde{\lambda},\tilde{\mu}}
\mathbf{u}^{\omega_0}}-\mathcal{L}^{\tilde{\lambda},\tilde{\mu}}
\mathbf{u}^{\omega}\cdot
\overline{\mathbf{u}^{\omega_0}}\\
&=\int_{\p
D}\mathbf{u}^\omega\cdot\pd{\overline{\mathbf{u}^{\omega_0}}}{\tilde{\nu}}
=(\omega^2-\omega_0^2)\int_{\p D}\psi^\omega
\cdot\pd{\overline{\mathbf{u}^{\omega_0}}}{\tilde{\nu}}.
\end{align*}
Dividing by $\omega^2-\omega_0^2$ and letting $\omega\rightarrow
\omega_0$, we have
\begin{align*}
\int_D |\mathbf{u}^{\omega_0}|^2 =\int_{\p D}\psi^{\omega_0}
\cdot\pd{\overline{\mathbf{u}^{\omega_0}}}{\tilde{\nu}}.
\end{align*}
Therefore we conclude that $\psi^{\omega_0}$ is not identically
zero. This completes the proof. \qed

By Lemma \ref{invert} and the fact that $\tilde{\Scal}^{\omega}$
is Fredholm, we know that $\Acal^{\alpha,\omega}_0$ is normal.
Moreover, Lemma \ref{rankone} says that the multiplicity of
$\Acal^{\alpha,\omega}_0$ at each eigenvalue $\omega_0^2$ of
$-\mathcal{L}^{\tilde{\lambda},\tilde{\mu}}$ is equal to the
dimension of $\ker\tilde{\Scal}^{\omega_0}$. Combining this fact
with Theorem \ref{rouche}, we obtain the following lemma.
\begin{lem}\label{eigen-neumann8}
For each eigenvalue $\omega_0^2$ of
$-\mathcal{L}^{\tilde{\lambda},\tilde{\mu}}$ and sufficiently
large $\mu$, there exists a small neighborhood $V$ of $\omega_0>0$
such that $\mathcal{A}^{\alpha,\omega}$ is normal with respect to
$\p V$ and $\mathcal{M}(\mathcal{A}^{\alpha,\omega}, \p V)=\dim
\ker\tilde{\Scal}^{\omega_0}$.
\end{lem}

Let $\omega_0^2$ (with $\omega_0 >0$) be a simple eigenvalue of
$-\mathcal{L}^{\tilde{\lambda},\tilde{\mu}}$ in $D$ with the
Dirichlet boundary condition. There exists a unique eigenvalue
$(\omega^\alpha_\mu)^2$ (with $\omega^\alpha_\mu>0$) of
\eqnref{transprime3} lying in a small complex neighborhood $V$ of
$\omega_0$. Combining the generalized Rouch{\'e}'s theorem  with
Lemma \ref{asympta} we are now able to derive complete asymptotic
formulae for the characteristic values of $\omega \mapsto
\Acal^{\alpha,\omega}$. Theorem \ref{formula} yields that
 \beq
 \omega^\alpha_\mu - \omega_0 = \frac{1}{2\sqrt{-1}\pi} \mbox{ {\rm tr} }\int_{\partial
 V} (\omega- \omega_0)(\Acal^{\alpha,\omega})^{-1}\frac{d}{d\omega}
 \Acal^{\alpha,\omega}d\omega.
 \eeq
Then we obtain the following complete asymptotic expansion for the
eigenvalue perturbations $\omega_\mu^\alpha - \omega_0$.

\begin{thm} \label{thmd3eigenalpha}
Suppose $\alpha \neq 0$. Then the following asymptotic expansion
holds:
\begin{align}\label{mainexp2}
\omega^\alpha_\mu-\omega_0&=\frac{1}{2\pi
\sqrt{-1}}\sum_{p=1}^{+\infty}\frac{1}{p}\sum_{n=p}^{+\infty}
\frac{1}{\mu^n} \tr\int_{\p V} \mathcal{B}^{\alpha,\omega}_{n,p}
d\omega ,
\end{align}
where
\begin{align}
\mathcal{B}^{\alpha,\omega}_{n,p}=(-1)^p\sum_{
n_1+\cdots+n_p=n\atop n_i\geq 1
}(\mathcal{A}^{\alpha,\omega}_0)^{-1}\mathcal{A}^{\alpha,\omega}_{n_1}\cdots(\mathcal{A}^{\alpha,\omega}_0)^{-1}\mathcal{A}^{\alpha,\omega}_{n_p}.
\end{align}
\end{thm}

\pf For sufficiently large $\mu$, the following Neumann series
converges uniformly with respect to the variable $\omega \in \p
V$:
\begin{align*}
(\mathcal{A}^{\alpha,\omega})^{-1}=\sum_{p=0}^{+\infty}
\bigg[(\mathcal{A}^{\alpha,\omega}_0)^{-1}
(\mathcal{A}^{\alpha,\omega}_0-\mathcal{A}^{\alpha,\omega})
\bigg]^p(\mathcal{A}^{\alpha,\omega}_0)^{-1}.
\end{align*}
By \eqnref{trace} and the relation
\begin{align*}
\frac{d}{d\omega}(\mathcal{A}^{\alpha,\omega}_0)^{-1}=-(\mathcal{A}^{\alpha,\omega}_0)^{-1}\frac{d}{d\omega}\mathcal{A}^{\alpha,\omega}_0(\mathcal{A}^{\alpha,\omega}_0)^{-1},
\end{align*}
we get
\begin{align*}
&\frac{1}{2\pi \sqrt{-1}}\tr\int_{\p
V}(\omega-\omega_0)\frac{1}{p}\frac{d}{d\omega}
\bigg[(\mathcal{A}^{\alpha,\omega}_0)^{-1}
(\mathcal{A}^{\alpha,\omega}_0-\mathcal{A}^{\alpha,\omega}) \bigg]^p d\omega\\
&=\frac{1}{2\pi \sqrt{-1}}\tr\int_{\p
V}(\omega-\omega_0)\bigg[(\mathcal{A}^{\alpha,\omega}_0)^{-1}
(\mathcal{A}^{\alpha,\omega}_0-\mathcal{A}^{\alpha,\omega})\bigg]^{p-1}(\mathcal{A}^{\alpha,\omega}_0)^{-1}\frac{d}{d\omega}(\mathcal{A}^{\alpha,\omega}_0-\mathcal{A}^{\alpha,\omega})d\omega\\
&~~-\frac{1}{2\pi \sqrt{-1}}\tr\int_{\p
V}(\omega-\omega_0)\bigg[(\mathcal{A}^{\alpha,\omega}_0)^{-1}
(\mathcal{A}^{\alpha,\omega}_0-\mathcal{A}^{\alpha,\omega})\bigg]^{p}(\mathcal{A}^{\alpha,\omega}_0)^{-1}\frac{d}{d\omega}\mathcal{A}^{\alpha,\omega}_0
d\omega.
\end{align*}
Summing over $p$, we obtain
\begin{align*}
&\frac{1}{2\pi \sqrt{-1}}\sum^{+\infty}_{p=1}\tr\int_{\p
V}(\omega-\omega_0)\frac{1}{p}\frac{d}{d\omega}\bigg[(\mathcal{A}^{\alpha,\omega}_0)^{-1}
(\mathcal{A}^{\alpha,\omega}_0-\mathcal{A}^{\alpha,\omega})\bigg]^p d\omega\\
&=-\frac{1}{2\pi \sqrt{-1}}\sum^{+\infty}_{p=0}\tr\int_{\p
V}(\omega-\omega_0)\bigg[(\mathcal{A}^{\alpha,\omega}_0)^{-1}
(\mathcal{A}^{\alpha,\omega}_0-\mathcal{A}^{\alpha,\omega})\bigg]^p(\mathcal{A}^{\alpha,\omega}_0)^{-1}\frac{d}{d\omega}\mathcal{A}^{\alpha,\omega}
d\omega\\
&~~+\frac{1}{2\pi \sqrt{-1}}\tr\int_{\p
V}(\omega-\omega_0)(\mathcal{A}^{\alpha,\omega}_0)^{-1}\frac{d}{d\omega}\mathcal{A}^{\alpha,\omega}_0
d\omega.
\end{align*}
Since \begin{align*} \frac{1}{2\pi \sqrt{-1}}\tr\int_{\p
V}(\omega-\omega_0)(\mathcal{A}^{\alpha,\omega}_0)^{-1}\frac{d}{d\omega}\mathcal{A}^{\alpha,\omega}_0
d\omega=0,
\end{align*}
and
\begin{align*}
&\frac{1}{2\pi \sqrt{-1}}\tr\int_{\p
V}(\omega-\omega_0)\frac{d}{d\omega}\bigg[(\mathcal{A}^{\alpha,\omega}_0)^{-1}
(\mathcal{A}^{\alpha,\omega}_0-\mathcal{A}^{\alpha,\omega})\bigg]^p
d\omega\\
&=-\frac{1}{2\pi \sqrt{-1}}\tr\int_{\p
V}\bigg[(\mathcal{A}^{\alpha,\omega}_0)^{-1}
(\mathcal{A}^{\alpha,\omega}_0-\mathcal{A}^{\alpha,\omega})\bigg]^p
d\omega,
\end{align*}
we have
\begin{align*}
&\frac{1}{2\pi \sqrt{-1}}\sum^{+\infty}_{p=1}\tr\int_{\p
V}\frac{1}{p}\bigg[(\mathcal{A}^{\alpha,\omega}_0)^{-1}
(\mathcal{A}^{\alpha,\omega}_0-\mathcal{A}^{\alpha,\omega})\bigg]^p
d\omega\\
&=\frac{1}{2\pi \sqrt{-1}}\sum^{+\infty}_{p=0}\tr\int_{\p
V}(\omega-\omega_0)\bigg[(\mathcal{A}^{\alpha,\omega}_0)^{-1}
(\mathcal{A}^{\alpha,\omega}_0-\mathcal{A}^{\alpha,\omega})\bigg]^p(\mathcal{A}^{\alpha,\omega}_0)^{-1}\frac{d}{d\omega}\mathcal{A}^{\alpha,\omega}
d\omega.\\
&=\frac{1}{2\pi \sqrt{-1}}\tr\int_{\p
V}(\omega-\omega_0)(\mathcal{A}^{\alpha,\omega})^{-1}\frac{d}{d\omega}\mathcal{A}^{\alpha,\omega}
d\omega.
\end{align*}
By expanding $\big[(\mathcal{A}^{\alpha,\omega}_0)^{-1}
(\mathcal{A}^{\alpha,\omega}_0-\mathcal{A}^{\alpha,\omega})\big]^p$,
we obtain the desired result. \qed

\subsubsection{{The case $\alpha =0$.}} We now deal with the
periodic case ($\alpha =0$). By \eqnref{gexp2} we have
\begin{equation}
\Acal^{0,\omega} = \Acal^{0,\omega}_0 +\sum_{l=1}^{+\infty}
\frac{1}{\mu^l}\Acal^{0,\omega}_l,
\end{equation}
where
 \beq \label{zeroomzero}
 \Acal^{0,\omega}_0= \left(\begin{array}{cc}
 \tilde{\Scal}^{\omega}  & -\ds\frac{1}{\omega^2}\ds\int_{\p
 D}\cdot~
 d\sigma
 \\ \nm \ds \frac{1}{2}I- (\tilde{\Kcal}^{\omega})^* & \ds
\frac{1}{2}I+ (\Kcal^{0,0})^*
 \end{array} \right),
 \eeq
and, for $l\geq 1$,
 \beq \label{zeroomell}
 \Acal^{0,\omega}_l=
 \left(\begin{array}{cc} 0 & - \Scal_{l}^{0,\omega} \\
 \nm  0 & \frac{1}{\mu}(\Kcal_{l+1}^{0,\omega})^*
 \end{array} \right).
 \eeq

Here we consider the following eigenvalue problem
\begin{align}\label{neweq}
\begin{cases} \ds
(\mathcal{L}^{\tilde{\lambda},\tilde{\mu}}+\omega^2)\mathbf{u} =0
\quad &
\mbox{in }~ D, \\
\nm \mathbf{u}+\ds\frac{1}{|Y\setminus D|}\int_D \mathbf{u}=0
\quad &\mbox{on } \p D.
\end{cases}
\end{align}
We note that it has a discrete spectrum and its eigenvalues are
nonnegative since we have
\begin{align*}
\int_D
\tilde{\lambda}|\nabla\cdot\mathbf{u}|^2+\frac{\tilde{\mu}}{2}|\nabla\mathbf{u}+\nabla\mathbf{u}^t|^2&=\int_{\p
D}\mathbf{u}\cdot\pd{\overline{\mathbf{u}}}{\tilde{\nu}}-\int_D
\mathbf{u}\cdot\mathcal{L}^{\tilde{\lambda},\tilde{\mu}}\overline{\mathbf{u}}\\
&=-\frac{1}{|Y\setminus D|}\int_D \mathbf{u}\cdot\int_{\p
D}\pd{\overline{\mathbf{u}}}{\tilde{\nu}}+\overline{\omega}^2\int_D
|\mathbf{u}|^2\\
&=\frac{\overline{\omega}^2}{|Y\setminus D|}\left|\int_D
\mathbf{u}~\right|^2+ \overline{\omega}^2\int_D |\mathbf{u}|^2.
\end{align*}
The eigenvalue of \eqnref{neweq}  is related with the
characteristic value of $\Acal^{0,\omega}$ as follows.
 \begin{lem} \label{lemnew}
The equation \eqnref{neweq} has a nonzero solution  if and only if
$\omega$ is a characteristic value of the operator-valued function
$\Acal^{0,\omega}_0$.
\end{lem}
\pf Suppose that there exists a nonzero pair $(\phi, \psi)$ such
that
 $$
 \Acal^{0,\omega}_0 \begin{pmatrix} \phi \\ \psi
 \end{pmatrix}=0,
 $$
or equivalently
 \begin{align}
 \tilde{\Scal}^{\omega} \phi - \frac{1}{\omega^2} \int_{\p
 D} \psi d\sigma &=0  \quad \mbox{on } \p D, \label{somegad1} \\
 \left( \frac{1}{2}I- (\tilde{\Kcal}^{\omega})^* \right) \phi +
 \left(\frac{1}{2}I+ (\Kcal^{0,0})^* \right) \psi &=0 \quad \mbox{on } \p D.
 \label{somega2}
 \end{align}
In particular,  $\phi$ is nonzero by the invertibility of
$\frac{1}{2}I+ (\Kcal^{0,0})^*$. Let
$\mathbf{u}:=\tilde{\Scal}^{\omega}\phi$. Then we have
\begin{align*}
\ds\frac{1}{|Y\setminus D|}\int_D
\mathbf{u}&=-\ds\frac{1}{\omega^2|Y\setminus D|}\int_{\p D}
\pd{\mathbf{u}}{\tilde{\nu}}\\
&=-\ds\frac{1}{\omega^2|Y\setminus D|}\int_{\p
D}\big(-\frac{1}{2}I+(\tilde{\Kcal}^{\omega})^*\big)\phi\\
&=-\ds\frac{1}{\omega^2|Y\setminus D|}\int_{\p D}
\left(\frac{1}{2}I+ (\Kcal^{0,0})^*\right)\psi\\
&=-\ds\frac{1}{\omega^2}\int_{\p D}\psi,
\end{align*}
where the last equality follows from \eqnref{d002}. Therefore by \eqnref{somegad1} $\mathbf{u}$ is a nonzero solution
to \eqnref{neweq}.

Suppose that \eqnref{neweq} has nonzero solution $\mathbf{u}$. Following the same argument as in the proof of Theorem \ref{repthm},
we can see that there exists
$\phi$ such that
\begin{equation}\label{phipsi}
 \begin{cases}
\tilde{\Scal}^{\omega}\phi=\mathbf{u}|_{\partial D}& \text{ on }\partial D,\\
\nm
\ds\big(-\frac{1}{2}I+(\tilde{\Kcal}^{\omega})^*\big)\phi=\frac{\partial
\mathbf{u}}{\partial\tilde{\nu}}& \text{ on }\partial D.
\end{cases}
\end{equation}
If we set \begin{align*} \psi=(\frac{1}{2}I+
(\Kcal^{0,0})^*)^{-1}\left(\pd{\mathbf{u}}{\tilde{\nu}}\right),
\end{align*}
then $(\phi,\psi)$ satisfies
$$ \Acal^{0,\omega}_0 \begin{pmatrix} \phi \\ \psi
\end{pmatrix}=0.$$
This completes the proof. \qed

 We also have the following Lemma.
\begin{lem}
Every eigenvector of $\Acal^{0,\omega}_0$ has rank one.
\end{lem}
\pf Suppose that $\begin{pmatrix} \phi \\ \psi
\end{pmatrix}$ is an eigenvector of $\Acal^{0,\omega}_0$ with rank
$m$ associated with characteristic value $\omega_0$, {\it i.e.},
there exist $\phi^\omega$ and $\psi^\omega$, holomorphic as
functions of $\omega$, such that $\phi^{\omega_0}=\phi$,
$\psi^{\omega_0}=\psi$, and
$$ \Acal^{0,\omega}_0 \begin{pmatrix} \phi^\omega \\ \psi^\omega
\end{pmatrix}=(\omega-\omega_0)^m\begin{pmatrix} \widetilde{\phi}^\omega \\
\widetilde{\psi}^\omega
\end{pmatrix}, $$
for some $\begin{pmatrix} \widetilde{\phi}^\omega \\
\widetilde{\psi}^\omega
\end{pmatrix} \in L^2(\p D)^2$. In other words, the following identities hold on $\p D$:
 \begin{align*}
 \tilde{\Scal}^{\omega} \phi^\omega - \frac{1}{\omega^2} \int_{\p
 D} \psi^\omega d\sigma &= (\omega -\omega_0)^m \widetilde{\phi}^\omega , \\
 \left( \frac{1}{2}I- (\tilde{\Kcal}^{\omega})^* \right) \phi^\omega +
 \left(\frac{1}{2}I+ (\Kcal^{0,0})^* \right) \psi^\omega &=(\omega -\omega_0)^m \widetilde{\psi}^\omega .
 \end{align*}
It then follows from \eqnref{d002} that
 \begin{align*}
 & \tilde{\Scal}^{\omega}\phi^\omega-\frac{1}{|Y\setminus D|\omega^2}
 \int_{\p D}\big(-\frac{1}{2}I+(\tilde{\Kcal}^{\omega})^*\big)\phi^\omega
 d\sigma \\
 & =\tilde{\Scal}^{\omega}\phi^\omega-\frac{1}{|Y\setminus D|\omega^2}
 \int_{\p D}\big( \frac{1}{2}I+ (\Kcal^{0,0})^*  \big)\psi^\omega
 d\sigma + \frac{(\omega -\omega_0)^m}{|Y\setminus D|\omega^2}
 \int_{\p D}\widetilde{\psi}^\omega d\sigma \\
 & =\tilde{\Scal}^{\omega}\phi^\omega-\frac{1}{\omega^2}
 \int_{\p D} \psi^\omega d\sigma + \frac{(\omega -\omega_0)^m}{|Y\setminus D|\omega^2}
 \int_{\p D}\widetilde{\psi}^\omega d\sigma \\
 &= (\omega-\omega_0)^m\left(\tilde{\phi}^\omega + \frac{1}{|Y\setminus
 D|\omega^2}\int_{\p D}\widetilde{\psi}^\omega d\sigma \right).
 \end{align*}
Let
 $$
 \eta^\omega: =\left(\tilde{\phi}^\omega + \frac{1}{|Y\setminus
 D|\omega^2} \int_{\p D}\widetilde{\psi}^\omega d\sigma \right) \quad \mbox{and}\quad
 \mathbf{u}^\omega:=\tilde{\Scal}^{\omega}\phi^\omega.
 $$
Then $\mathbf{u}^\omega$ satisfies
\begin{equation*}
\begin{cases}
\ds (\mathcal{L}^{\tilde{\lambda},\tilde{\mu}}+\omega^2)\mathbf{u}^\omega=0\quad\mbox{in}~D,\\
\nm \mathbf{u}^\omega=\ds\frac{1}{|Y\setminus D|\omega^2}\int_{\p
D}\pd{\mathbf{u}^\omega}{\tilde{\nu}}d\sigma+(\omega-\omega_0)^m\eta^\omega\quad\mbox{on}~\p
D.
\end{cases}
\end{equation*}
By Green's formula, we have
\begin{align*}
&(\omega^2-\omega_0^2)\int_D \mathbf{u}^\omega
\cdot\overline{\mathbf{u}^{\omega_0}}\\&=\int_{\p
D}\mathbf{u}^\omega\cdot\pd{\overline{\mathbf{u}^{\omega_0}}}{\tilde{\nu}}-\overline{\mathbf{u}^{\omega_0}}\cdot\pd{\mathbf{u}^\omega}{\tilde{\nu}}
d\sigma\\
&=\left(\frac{1}{\omega^2}-\frac{1}{\omega_0^2}\right)\frac{1}{|Y\setminus
D|}\int_{\p D}\pd{\mathbf{u}^\omega}{\tilde{\nu}}d\sigma
\cdot\int_{\p
D}\pd{\overline{\mathbf{u}^{\omega_0}}}{\tilde{\nu}}d\sigma+(\omega-\omega_0)^m\int_{\p
D}\eta^\omega
\cdot\pd{\overline{\mathbf{u}^{\omega_0}}}{\tilde{\nu}}d\sigma.
\end{align*}
Dividing by $\omega^2-\omega_0^2$ and letting $\omega\rightarrow
\omega_0$, we obtain
\begin{align*}
\int_D \big|\mathbf{u}^{\omega_0}\big|^2+\frac{1}{2|Y\setminus D|
\omega_0^4}\Bigg|\int_{\p
D}\pd{\mathbf{u}^{\omega_0}}{\tilde{\nu}}d\sigma\Bigg|^2=\lim_{\omega\rightarrow\omega_0}\frac{(\omega-\omega_0)^m}{\omega^2-\omega_0^2}\int_{\p
D}\eta^{\omega_0}\cdot\pd{\overline{\mathbf{u}^{\omega_0}}}{\tilde{\nu}}d\sigma.
\end{align*}
Since the term on the left is nonzero, we conclude that $m=1$.
This completes the proof. \qed

 Analogously to Theorem \ref{thmd3eigenalpha}, the following asymptotic
formula for $\alpha=0$ holds.
\begin{thm} \label{thmd3eigenalphazero}
Suppose $\alpha = 0$. Let $\widetilde{\omega}_0^2$ (with
$\widetilde{\omega}_0>0$) be a simple eigenvalue of
\eqnref{neweq}. Then there exists a unique characteristic value
$\omega^0_\mu>0$ of $\Acal^{0,\omega}$ lying in a small complex
neighborhood $V$ of $\widetilde{\omega}_0$ and the following
asymptotic expansion holds:
\begin{align} \label{formula2}
\omega^0_\mu-\tilde{\omega}_0&=\frac{1}{2\pi
\sqrt{-1}}\sum_{p=1}^{+\infty}\frac{1}{p}\sum_{n=p}^{+\infty}
\frac{1}{\mu^n} \tr\int_{\p V} \mathcal{B}^{\omega}_{n,p} d\omega,
\end{align}
where
\begin{align}
\mathcal{B}^{\omega}_{n,p}=(-1)^p\sum_{ n_1+\cdots+n_p=n\atop
n_i\geq 1
}(\mathcal{A}^{0,\omega}_0)^{-1}\mathcal{A}^{0,\omega}_{n_1}\cdots(\mathcal{A}^{0,\omega}_0)^{-1}\mathcal{A}^{0,\omega}_{n_p}.
\end{align}
\end{thm}

\subsubsection{{The case when $|\alpha|$ is of order of $1/\sqrt{\mu}$.}} In
this subsection we derive an asymptotic expansion  which is valid
for $|\alpha|$ of order $O(1/\sqrt{\mu})$, not just for fixed
$\alpha \neq 0$ or $\alpha =0$, as has been considered in the
previous subsections and give the limiting behavior of
$\omega^\alpha_\mu$ in this case.

Suppose that $|\alpha|^2 \mu$ goes to $0< \tau < +\infty$ and
denote $\bar{\bar{\tau}} =(\tau_{ij})_{i,j=1,2}$, with $\tau_{ij}
= \lim \alpha_i \alpha_j \mu$. Then, following the same arguments
as those in the proof of Lemma \ref{lemnew}, we can show that the
following problem
\begin{align}\label{neweqn} \begin{cases} \ds
(\mathcal{L}^{\tilde{\lambda},\tilde{\mu}}+\omega^2)\mathbf{u} =0
\quad &
\mbox{in }~ D, \\
\nm \mathbf{u}+\ds  \frac{1}{|Y\setminus D|} \bigg[\frac{1}{1 -
\frac{\tau}{\omega^2}} I + \frac{\bar{\bar{\tau}}}{ (\omega^2 - 2
\tau) (1 - \frac{\tau}{\omega^2})} \bigg] \int_D \mathbf{u}=0
\quad &\mbox{on } \p D,
\end{cases}
\end{align}
has a nontrivial solution if and only if $\omega$ is a real
characteristic value of the operator-valued function
$$
 \Acal^{\omega}_\tau = \left(\begin{array}{cc}
 \widetilde{\Scal}^{\omega}  & -\ds\frac{1}{\omega^2- \tau} \bigg[ I + \frac{\bar{\bar{\tau}}} { \omega^2 - 2
\tau} \bigg] \ds\int_{\p
 D}\cdot~
 d\sigma
 \\ \nm \ds \frac{1}{2}I- (\widetilde{\Kcal}^{\omega})^* & \ds
\frac{1}{2}I+ (\Kcal^{0,0})^*
 \end{array} \right).
 $$
We can also show that any eigenvector of $\Acal^{\omega}_\tau$ has
rank one. Then if we denote by $\widetilde{\omega}_\tau$ a simple
eigenvalue of \eqnref{neweqn} then
\beq \label{alphamu}
 \omega_\mu^\alpha - \widetilde{\omega}_\tau = \frac{1}{2\pi \sqrt{-1}} \frac{1}{\mu}
 \int_{\partial V} (\Acal^{\omega}_\tau)^{-1} \Acal^{0,\omega}_1\; d\omega
 + O(1/\mu^2),
\eeq
as $\mu \rightarrow +\infty$ and $|\alpha|^2 \mu
\rightarrow \tau$. Here $\Acal^{0,\omega}_1$ is defined by the
same formula as in the previous subsection.

Not surprisingly, this asymptotic expansion tends continuously to
\eqnref{mainexp2} and \eqnref{formula2} as $\tau$ goes to $+
\infty$ or $0$, respectively.

\subsection{Derivation of the leading order terms}

For $\alpha\ne0$, let us write down explicitly the leading order
term in the expansion of $\omega_\mu^\alpha - \omega_0$. We first observe that
\begin{equation} \label{aoo} (\Acal^{\alpha,\omega}_0)^{-1} =
\begin{pmatrix} (\widetilde{\Scal}^{\omega})^{-1}  & &  0
\\ \nm \ds
\left(\frac{1}{2}I+ (\Kcal^{\alpha,0})^*\right)^{-1}\left(
\frac{1}{2}I- (\widetilde{\Kcal}^{\omega})^*
\right)(\widetilde{\Scal}^{\omega})^{-1} & & \ds \left(\frac{1}{2}I+
(\Kcal^{\alpha,0})^*\right)^{-1}
\end{pmatrix}. \end{equation}

Next, we prove the following Lemma.
\begin{lem}
Let $\mathbf{u}_0$ be an eigenvector associated to the simple
eigenvalue $\omega_0^2$ and let $\vp:=
\pd{\mathbf{u}_0}{\tilde{\nu}}|_{-}$ on $\p D$. Then we have, in a
neighborhood of $\omega_0$,
\beq \label{334}
(\widetilde{\Scal}^{\omega})^{-1}= \frac{1}{\omega-\omega_0} T+
\Qcal^{\omega}
\eeq
where $\Qcal^{\omega}$ is operators in $\mathcal{L}(H^2(\p D)^2, L^2(\p D)^2)$ holomorphic in $\omega$, and
$T$ is defined by
\beq \label{Tdef}
T(f):=-\frac{\langle f,\vp \rangle \vp}{2\omega_0\ds\int_D|\mathbf{u}_0|^2},
\eeq
where $\langle \ , \ \rangle$ is the inner product on $L^2(\p D)^2$.
\end{lem}
\pf By Lemma \ref{rankone}, there are operators $T$ and $Q^\omega$ in $\mathcal{L}(H^2(\p D)^2, L^2(\p D)^2)$
such that $(\widetilde{\Scal}_D^{\omega})^{-1}$ takes the form
\beq \label{336}
(\widetilde{\Scal}^{\omega})^{-1}= \frac{1}{\omega-\omega_0} T
 + \Qcal^{\omega}
\eeq
where $\Qcal^{\omega}$ is holomorphic in $\omega$. Since
\beq
\mbox{Id} = (\widetilde{\Scal}^{\omega})(\widetilde{\Scal}^{\omega})^{-1} = \frac{1}{\omega-\omega_0} \widetilde{\Scal}^{\omega} T
+ \widetilde{\Scal}^{\omega} \Qcal^{\omega},
\eeq
by letting $\omega \to \omega_0$, we have
\beq \label{sazero}
\widetilde{\Scal}_D^{\omega_0} T=0.
\eeq
Similarly, we can show that
\beq \label{aszero}
T \widetilde{\Scal}_D^{\omega_0}=0.
\eeq
It then follows from \eqnref{sazero} and \eqnref{aszero} that
$\im A=\ker\widetilde{\Scal}^{\omega_0}= \mbox{ span} \{
\varphi \}$ and $\ker A= \im \widetilde{\Scal}^{\omega_0}= \mbox{
span} \{ \varphi \}^\perp$. Here $\mbox{ span} \{\varphi \}$
denotes the vector space spanned by $\varphi$. Therefore
\beq \label{TClan}
T=C \langle\cdot, \vp \rangle \vp
\eeq
for some constant $C$.

By Green's formula, we have for $x \in D$
\begin{align}
\widetilde{\Scal}^\omega(\vp)(x) & = \widetilde{\Scal}^\omega \left( \pd{\mathbf{u}_0}{\nu} \Big|_{-} \right)(x) -
\widetilde{\Dcal}^\omega \left( \mathbf{u}_0 \right)(x) \nonumber \\
&=(\omega^2-\omega_0^2)\int_D
\widetilde{\Gamma}^\omega(x-y)\mathbf{u}_0(y)dy - \mathbf{u}_0(x). \label{341}
\end{align}
In particular, we get
\beq \label{omeome1}
\widetilde{\Scal}^\omega(\vp)(x) = (\omega^2-\omega_0^2)\int_D
\widetilde{\Gamma}^\omega(x-y)\mathbf{u}_0(y)dy, \quad x \in \p D.
\eeq
By expanding $\widetilde{\Gamma}^\omega(x-y)$ in $\omega$, we now have
\beq \label{omeome2}
\widetilde{\Scal}^\omega(\vp)(x) =2\omega_0(\omega-\omega_0)\int_D
\widetilde{\Gamma}^{\omega_0}(x-y)\mathbf{u}_0(y)dy+(\omega-\omega_0)^2A^\omega,
\eeq
for some function $A^\omega$ holomorphic in $\omega$. Therefore, we have
\beq \label{twoome1}
(\widetilde{\Scal}^\omega)^{-1}\left(2\omega_0\int_D
\widetilde{\Gamma}^{\omega_0}(x-y)\mathbf{u}_0(y)dy\right)=\frac{1}{\omega-\omega_0}\vp
+ B^\omega,
\eeq
where $B^\omega$ is holomorphic in $\omega$, which together with \eqnref{336} implies that
\beq \label{twoome2}
T\left(2\omega_0\int_D
\widetilde{\Gamma}^{\omega_0}(\cdot-y)\mathbf{u}_0(y)dy\right)=\vp.
\eeq
Note that if we take $\omega=\omega_0$ in \eqnref{341}, then
 \beq \label{341-1}
 \mathbf{u}_0(x)= -\widetilde{\Scal}^{ \omega_0}(\vp)(x), \quad x \in D.
 \eeq
It then follows from \eqnref{TClan} and \eqnref{twoome2} that
\begin{align*}
1 & =C\left\langle 2\omega_0\int_D
\widetilde{\Gamma}^{\omega_0}(x-y)\mathbf{u}_0(y)dy,\vp\right\rangle \\
&= 2C \omega_0 \left\langle \mathbf{u}_0, \widetilde{\Scal}^\omega \vp\right\rangle=-2C\omega_0\int_D|\mathbf{u}_0|^2.
\end{align*}
This completes the proof. \qed

Because of \eqnref{341-1}, we have
\beq \label{347}
(\frac{1}{2}I- (\widetilde{\Kcal}^{\omega_0})^*) (\vp) = \vp \quad \mbox{on } \p D.
\eeq
Observe from \eqnref{317} and \eqnref{318} that the diagonal elements of $(\Acal^{\alpha,\omega}_{0})^{-1}
\Acal^{\alpha,\omega}_{1}$ are $0$ and
\beq
- \big( \frac{1}{2}I+ (\Kcal^{\alpha,0})^* \big)^{-1}
\big(\frac{1}{2}I- (\widetilde{\Kcal}^{ \omega_0})^* \big) (\widetilde \Scal^\omega)^{-1} \Scal_1^{\alpha, \omega}
+ \mbox{an operator holomorphic in } \omega.
\eeq
The identity \eqnref{314} implies that $\Scal_1^{\alpha, \omega}= \mu \Scal^{\alpha, 0}$, and hence it follows from
\eqnref{334} that
 $$
 \frac{1}{2\pi \sqrt{-1}} \mbox{ tr }
 \int_{\partial{V}} (\Acal^{\alpha,\omega}_{0})^{-1}
 \Acal^{\alpha,\omega}_{1} d\omega=  -\mu \mbox{tr} \Big[ T \Scal^{\alpha,0}
\big( \frac{1}{2}I+ (\Kcal^{\alpha,0})^* \big)^{-1}
\big(\frac{1}{2}I- (\widetilde{\Kcal}^{ \omega_0})^* \big) \Big].
 $$
Since $\im T=\mbox{span}\{ \vp \}$, it follows from \eqnref{347} that
\begin{align}
& \mbox{tr} \Big[ T \Scal^{\alpha,0}
\big( \frac{1}{2}I+ (\Kcal^{\alpha,0})^* \big)^{-1}
\big(\frac{1}{2}I- (\widetilde{\Kcal}^{ \omega_0})^* \big) \Big] \nonumber \\
&= \frac{\left\langle \Big[ T \Scal^{\alpha,0}
\big( \frac{1}{2}I+ (\Kcal^{\alpha,0})^* \big)^{-1}
\big(\frac{1}{2}I- (\widetilde{\Kcal}^{ \omega_0})^* \big) \Big] \vp, \vp \right\rangle}{||\vp ||^2_{L^2(\p D)}} \nonumber \\
&= \frac{\left\langle \Big[ T \Scal^{\alpha,0}
\big( \frac{1}{2}I+ (\Kcal^{\alpha,0})^* \big)^{-1}
 \Big] \vp, \vp \right\rangle}{||\vp ||^2_{L^2(\p D)}} .
\end{align}
We set
 \begin{equation}
 \mathbf{v}_0(x):= \mu \Scal^{\alpha,0}
\left( \frac{1}{2}I+ (\Kcal^{\alpha,0})^* \right)^{-1} (\vp)(x),
\quad x \in Y\setminus \overline{D}.
 \end{equation}
Then $\mathbf{v}_0$ is the unique $\alpha$-quasi-periodic solution
to
 \begin{equation} \left\{
\begin{array}{ll}
\ds \mathcal{L}^{\lambda,\mu}\mathbf{v}_0 =0 \quad & \mbox{in }
Y\setminus
\overline{D}, \\
\nm \ds \pd{\mathbf{v}_0}{\nu} \Big|_{+} = \mu
\pd{\mathbf{u}_0}{\widetilde{\nu}} \Big|_{-} \quad &\mbox{on } \p
D,
\end{array}
\right. \end{equation} and
 \begin{align*}  \frac{1}{2\pi \sqrt{-1}} \tr
 \int_{\partial{V}} (\Acal^{\alpha,\omega}_{0})^{-1}
 \Acal^{\alpha,\omega}_{1} d\omega &= \frac{1}{||\vp ||^2_{L^2}} \langle \vp, T\mathbf{v}_0 \rangle \\
 &=\frac{1}{\mu} \frac{\ds \int_{Y\setminus \overline D}\lambda |\nabla\cdot\mathbf{v}_0|^2
 +\frac{\mu}{2}|\nabla\mathbf{v}_0 +\nabla\mathbf{v}_0^t|^2}{2\omega_0 \ds\int_D |\mathbf{u}_0|^2 } .
 \end{align*}
Thus the following corollary holds.
\begin{cor} Suppose $\alpha\ne 0$. Then the following asymptotic formula holds:
\beq
\omega_\mu^\alpha - \omega_0= \ds - \frac{1}{\mu}
\frac{\ds \int_{Y\setminus \overline D} \frac{\lambda}{\mu} |\nabla\cdot\mathbf{v}_0|^2
 +\frac{1}{2}|\nabla\mathbf{v}_0 +\nabla\mathbf{v}_0^t|^2}{2\omega_0 \ds\int_D |\mathbf{u}_0|^2 }
+ O\left(\frac{1}{\mu^2}\right) \quad \mbox{as } \mu \to + \infty.
\eeq
\end{cor}

When $\alpha =0$, it does not seem to be likely that we can explicitly compute the leading order term in a closed as in the
case $\alpha \neq 0$. But, let us now briefly explain how to compute the leading order term in
the asymptotic expansion of $\omega_\mu^0 - \widetilde{\omega}_0$.

Let $\mathbf{u}_0$ be the (normalized) eigenvector of
\eqnref{neweq} associated with the simple eigenvalue
$\widetilde{\omega}_0$. Let $(\widetilde{\phi}_0, \widetilde{\psi}_0)$ satisfy
\eqnref{phipsi} with $\mathbf{u}$ replaced by $\mathbf{u}_0$ and
$\omega = \widetilde{\omega}_0$. Since $\widetilde{\omega}_0$ is the only
simple pole in $V$ of the mapping $\omega \mapsto
(\Acal^{0,\omega}_0)^{-1}$, it can be proved that
\begin{align}
(\Acal^{0,\omega}_0)^{-1} &= \ds \frac{1}{\omega - \widetilde{\omega}_0}
\bigg(\frac{d}{d\omega} \Acal^{0,\omega}_0 \bigg|_{\omega =
\widetilde{\omega}_0}
\left( \begin{array}{l} \widetilde{\phi}_0 \\
\widetilde{\psi}_0 \end{array} \right) \cdot \left( \begin{array}{l} \widetilde{\phi}_0 \\
\widetilde{\psi}_0 \end{array} \right)  \bigg)^{-1}
\begin{pmatrix} \langle\cdot, \widetilde{\phi}_0 \rangle \; \widetilde{\phi}_0 & 0 \\ 0 &
\langle\cdot, \widetilde{\psi}_0 \rangle \; \widetilde{\psi}_0\end{pmatrix} \\
& \quad + \mbox{operator-valued function holomorphic in } \omega, \nonumber
\end{align}
which allows us to explicit the leading-order term in the
expansion of $\omega_\mu^0 - \widetilde{\omega}_0$. Similar
calculations and expressions in the transition region ($|\alpha| =
O(1/\sqrt{\mu})$) can be derived as well.

\subsection{Criterion for gap opening}

Following Hempel and Lienau \cite{hempel2}, we provide in this
subsection a criterion for gap opening in the spectrum of the
operator given by \eqnref{periodicop} as $\mu \rightarrow +
\infty$.

Let $\omega_j$ be the eigenvalues of
$-\mathcal{L}^{\tilde{\lambda},\tilde{\mu}}$ in $D$ with the
Dirichlet boundary condition. Let $\widetilde{\omega}_j$ denote
the eigenvalues of \eqnref{neweq}. We first prove the following
min-max characterization of $\omega_j$ and $\widetilde{\omega}_j$.

\begin{lem} \label{lemminmax}
The following min-max characterizations of $\omega^2_j$ and
$\widetilde{\omega}^2_j$ hold: \beq \label{minmax1} \omega^2_j =
\min_{N_j} \max_{ \mathbf{u} \in N_j, || \mathbf{u} || =1}
\widetilde{\mathbf{E}}(\mathbf{u},\mathbf{u}), \eeq and \beq
\label{minmax2} \widetilde{\omega}^2_j = \min_{N_j} \max_{
\mathbf{u} \in N_j, || \mathbf{u} || =1}
\frac{\widetilde{\mathbf{E}}(\mathbf{u},\mathbf{u})}{1 - |\int_D
\mathbf{u}|^2}, \eeq where the minimum is taken over all $j$
dimensional subspaces $N_j$ of $(H^1_0(D))^2$. Here $H^1_0(D)$ is
the set of all functions in $H^1(D)$ whose trace at $\partial D$
is zero and $\widetilde{\mathbf{E}}$ is given by \eqnref{E} with
$(\lambda, \mu)$ replaced by $(\tlambda, \tmu)$.
\end{lem}
\pf The identity \eqnref{minmax1} is well known. Note that if $\mathbf{v}$
satisfies the Dirichlet condition on $\partial D$, then
$\widetilde{\mathbf{v}}:= \mathbf{v} - \int_D \mathbf{v}$ satisfies
the boundary condition
\begin{equation} \label{pbc}
\widetilde{\mathbf{v}} + \frac{1}{|Y\setminus D|} \int_D \widetilde{\mathbf{v}} =0\quad
\mbox{on } \partial D.
\end{equation}
Conversely, if $\widetilde{\mathbf{v}}$ satisfies \eqnref{pbc}, then $\mathbf{v}:
= \widetilde{\mathbf{v}} + \frac{1}{|Y\setminus D|} \int_D
\widetilde{\mathbf{v}}$ obviously satisfies the Dirichlet boundary
condition.

Observe that the operator with the boundary condition in \eqnref{neweq} is not self-adjoint, and hence Poincare's min-max principle
can not be applied. So we now introduce an eigenvalue problem whose eigenvalues are exactly those of \eqnref{neweq}.
Let $\mathcal{H}=\mbox{span}\{H^2_0(D), 1_Y \}$ in $L^2(Y)$ where $H^2_0(D)$ is regarded as a subspace of $L^2(Y)$ by extending the functions
to be $0$ in $Y\setminus D$. Let $\mathcal{G}$
be the closure of $\mathcal{H}$ in $L^2(Y)$. Define the operator
$\mathbf{T} : \mathcal{H} \times \mathcal{H} \rightarrow \mathcal{G} \times \mathcal{G}$ by
\begin{align*}
\mathbf{T}\mathbf{u}=\begin{cases}
-\mathcal{L}^{\tilde{\lambda},\tilde{\mu}}\mathbf{u} \quad
&\mbox{on}
~ D,\\
\frac{1}{|Y\setminus D|}\int_D
\mathcal{L}^{\tilde{\lambda},\tilde{\mu}}\mathbf{u} \quad
&\mbox{on} ~Y\setminus D.
\end{cases}
\end{align*}
The constant value of $\mathbf{T}\mathbf{u}$ in $Y \setminus D$ was chosen so that $\int_Y \mathbf{T}\mathbf{u}=0$.
Then one can easily see that $\mathbf{T}$ is a densely defined self-adjoint operator on $\mathcal{H} \times \mathcal{H}$ and
\beq
\langle
\mathbf{T}\mathbf{u},\mathbf{v}\rangle_Y=\widetilde{\mathbf{E}}(\mathbf{u},\mathbf{v})\quad
\mbox{for} ~\mathbf{u},\mathbf{v}\in \mathcal{H} \times \mathcal{H}.
\eeq
One can also show that nonzero eigenvalues of $\mathbf{T}$
are eigenvalues of \eqnref{neweq}, and vice versa.

Let $M_j$ be a $j$-dimensional subspace of $\mathcal{H} \times \mathcal{H}$ perpendicular to
constant vectors which is eigenvectors corresponding to eigenvalue
zero. Then by Poincare's min-max principle we have
\begin{align*}
\widetilde{\omega}^2_j& = \min_{M_j} \max_{ \mathbf{u} \in M_j}
\frac{\langle \mathbf{T}\mathbf{u},\mathbf{u}\rangle_Y}{\langle
\mathbf{u},\mathbf{u}\rangle_Y}\\
&=\min_{M_j} \max_{ \mathbf{u} \in M_j}
\frac{\widetilde{\mathbf{E}}(\mathbf{u},\mathbf{u})}{\langle
\mathbf{u},\mathbf{u}\rangle_Y}\\
&=\min_{N_j} \max_{ \mathbf{v} \in N_j}
\frac{\widetilde{\mathbf{E}}(\mathbf{v}-\int_D
\mathbf{v},\mathbf{v}-\int_D \mathbf{v})}{\langle
\mathbf{v}-\int_D
\mathbf{v},\mathbf{v}-\int_D \mathbf{v}\rangle_Y}\\
&=\min_{N_j} \max_{ \mathbf{v} \in N_j}
\frac{\widetilde{\mathbf{E}}(\mathbf{v},\mathbf{v})}{\langle
\mathbf{v},\mathbf{v}\rangle_D-|\int_D \mathbf{v}|^2}.
\end{align*}
This complets the proof of the lemma. \qed

\begin{lem} \label{lemint}
The eigenvalues $\omega_j$ and $\widetilde{\omega}_j$ interlace in the following way:
\beq \label{358}
\omega_j \leq \widetilde{\omega}_j  \leq \omega_{j+2}, \quad j =1,2,\ldots.
\eeq
\end{lem}
\pf Lemma \ref{lemminmax} shows that the first inequality in
\eqnref{358} is trivial and we only have to prove the second one.
Let $\mathbf{u}_j$ denotes the normalized eigenvectors associated
with ${\omega}_j$. Let $N_{j+2}$ denotes the span of the
eigenvectors $\mathbf{u}_1, \ldots, \mathbf{u}_{j+2}$ and
$\widetilde{N}$ be the subspace of $N_{j+2}$ composed of all the
elements in $N_{j+2}$ which have zero integral over $D$. Since the
set of constant vectors has dimension two, $\widetilde{N}$ is of
dimension greater than $j$. Therefore, we have
$\widetilde{\omega}_j \leq \omega_{j+2}$, as desired. \qed

Since $0$ is eigenvalue of the periodic problem with multiplicity
2, formulae  \eqnref{mainexp2}, \eqnref{formula2}, and
\eqnref{alphamu} show that the spectral bands converge, as
$\mu\rightarrow\infty$, to \beq [0,\omega_1]\bigcup
[0,\omega_2]\bigcup_{j\geq1}[\widetilde{\omega}_j,\omega_{j+2}],
\eeq and hence we have a band-gap if and only if the following
holds:
\begin{equation} \label{crite1}
\omega_{j+1}<\widetilde{\omega}_j\quad \mbox{for some } j \quad (\mbox{\bf Criterion for gap opening}).
\end{equation}
Observe that by \eqnref{minmax1} and \eqnref{minmax2} the gap
opening criterion is equivalent to \beq \label{crite2}
\min_{N_{j+1}} \max_{ \mathbf{u} \in N_{j+1}, || \mathbf{u} || =1}
\widetilde{\mathbf{E}}(\mathbf{u},\mathbf{u}) < \min_{N_j} \max_{
\mathbf{u} \in N_j, || \mathbf{u} || =1}
\frac{\widetilde{\mathbf{E}}(\mathbf{u},\mathbf{u})}{1 - |\int_D
\mathbf{u}|^2}, \eeq where $N_j$ is a $j$-dimensional subspace of
$H_0^1(\p D)^2$.

It is not likely to find conditions on the inclusion $D$ so that
the gap-opening criterion is satisfied by rigorous analysis.
However, finding such conditions by all means such as numerical
computations will be of great importance. It should be emphasized
that the criterion \eqnref{crite1} is for the case when the matrix
and the inclusion have the same density, which we assumed to be
equal to  $1$.

\section{Gap opening criterion when densities are different}

We now consider periodic elastic composites such that the matrix
and the inclusion have different densities.

Suppose that the density of the matrix is $\rho$ while that of the
inclusion is $1$ (after normalization). The Lam\'e parameters are
the same as before. In this case, the first equation of the eigenvalue problem \eqnref{transprime3} is changed to
\beq
\mathcal{L}^{\lambda,\mu}\mathbf{u}+ \rho\omega^2\mathbf{u}=0, \quad\mbox{in } Y\setminus
\overline{D}.
\eeq
Hence we can show by exactly the
analysis that the asymptotic expansions \eqnref{mainexp2}, \eqnref{formula2}, and
\eqnref{alphamu} hold if we replace the operators \eqnref{zeroomzero} and \eqnref{318}
(and \eqnref{zeroomell}) with new operators (depending on the density $\rho$)
 \beq \label{zeroomzero2}
 \Acal^{0,\omega}_0= \left(\begin{array}{cc}
 \tilde{\Scal}^{\omega}  & -\ds\frac{1}{\rho\omega^2}\ds\int_{\p
 D}\cdot~
 d\sigma
 \\ \nm \ds \frac{1}{2}I- (\tilde{\Kcal}^{\omega})^* & \ds
\frac{1}{2}I+ (\Kcal^{0,0})^*
 \end{array} \right)
 \eeq
and
 \beq \label{newal}
 \Acal^{\alpha,\omega}_l= \rho^{l-1}
 \left(\begin{array}{cc} 0 & - \Scal_{l}^{\alpha,\omega} \\
 \nm  0 & \frac{\rho}{\mu}(\Kcal_{l+1}^{\alpha,\omega})^*
 \end{array} \right), \quad l \ge 1,
 \eeq
respectively, and the eigenvalue problems
\eqnref{neweq} and  \eqnref{neweqn} with the following eigenvalue problems
\begin{align}\label{neweq2} \begin{cases} \ds
(\mathcal{L}^{\tilde{\lambda},\tilde{\mu}}+\omega^2)\mathbf{u} =0
\quad &
\mbox{in }~ D, \\
\nm \mathbf{u}+\ds\frac{1}{\rho|Y\setminus D|}\int_D \mathbf{u}=0
\quad &\mbox{on } \p D.
\end{cases}
\end{align}
and
\begin{align}\label{neweqn2} \begin{cases} \ds
(\mathcal{L}^{\tilde{\lambda},\tilde{\mu}}+\omega^2)\mathbf{u} =0
\quad &
\mbox{in }~ D, \\
\nm \mathbf{u}+\ds  \frac{1}{\rho |Y\setminus D|} \bigg[\frac{1}{1 -
\frac{\tau}{\rho \omega^2}} I + \frac{\bar{\bar{\tau}}}{ (\rho\omega^2 - 2
\tau) (1 - \frac{\tau}{\rho\omega^2})} \bigg] \int_D \mathbf{u}=0
\quad &\mbox{on } \p D,
\end{cases}
\end{align}
respectively.

Let $\{ \widetilde{\omega}_j\}$ be the set of eigenvalues of
\eqnref{neweq2}. In order to express $\widetilde{\omega}_j$ using
the min-max principle, we define $\langle~ ,~\rangle_Y$  by \beq
\langle \mathbf{u} , \mathbf{v}\rangle_Y=\int_D \mathbf{u}\cdot
\mathbf{v }+\rho\int_{Y\setminus D} \mathbf{u}\cdot \mathbf{v}.
\eeq We also define, as before, $\mathbf{T}$ to be
\beq
\mathbf{T}\mathbf{u}=
\begin{cases}
\ds -\mathcal{L}^{\tilde{\lambda}, \tilde{\mu}} \mathbf{u}\quad &\mbox{on} ~D, \\
\nm
\ds \frac{1}{\rho|Y\setminus D|}\int_D \mathcal{L}^{\tilde{\lambda},
\tilde{\mu}} \mathbf{u} \quad &\mbox{on} ~Y\setminus D.
\end{cases}
\eeq
Then $\mathbf{T}$ is self-adjoint with respect to $\langle~
,~\rangle_Y$. By Poincare's min-max principle again, we have
\begin{align*}
\widetilde{\omega}_j^2&=\min_{M_j}\max_{\mathbf{u}\in M_j} \frac{\langle
\mathbf{T}\mathbf{u},
\mathbf{u}\rangle_Y}{\langle \mathbf{u}, \mathbf{u}\rangle_Y}\\
&=\min_{M_j}\max_{\mathbf{u}\in M_j} \frac{\widetilde{\mathbf{E}}(\mathbf{u},\mathbf{u})}{\langle \mathbf{u}, \mathbf{u}\rangle_Y}\\
&=\min_{N_j}\max_{\mathbf{v}\in N_j}
\frac{\widetilde{\mathbf{E}}(\mathbf{v},\mathbf{v})}{ \left\langle
\mathbf{v}- \frac{1}{|D|+\rho|Y\setminus D|}\int_D \mathbf{v},
\mathbf{v}-  \frac{1}{|D|+\rho|Y\setminus D|}\int_D
\mathbf{v} \right\rangle_Y}\\
&=\min_{N_j}\max_{\mathbf{v}\in N_j}
\frac{\widetilde{\mathbf{E}}(\mathbf{v},\mathbf{v})}{\langle
\mathbf{v},\mathbf{v}\rangle_D-  \frac{1}{|D|+\rho|Y\setminus
D|} \left| \int_D \mathbf{v} \right|^2},
\end{align*}
where $M_j$ and $N_j$ are the same as in the proof of Lemma
\ref{lemminmax}. Therefore, we have the following min-max
characterization of the eigenvalues of the problem
\eqnref{neweq2}:
 \beq
 \widetilde{\omega}_j^2=\min_{N_j}\max_{\mathbf{u}\in N_j}
 \frac{\widetilde{\mathbf{E}}(\mathbf{u},\mathbf{u})}{\langle
 \mathbf{u},\mathbf{u}\rangle_D-  \frac{1}{|D|+\rho|Y\setminus D|}
 \left| \int_D \mathbf{u} \right|^2}.
 \eeq
We then get a band-gap criterion for the different density case which is equivalent to \eqnref{crite1}:
 \beq \label{crite3}
 \min_{N_{j+1}} \max_{
 \mathbf{u} \in N_{j+1}, || \mathbf{u} || =1}
 \widetilde{\mathbf{E}}(\mathbf{u},\mathbf{u}) <
 \min_{N_j}\max_{\mathbf{u}\in N_j, || \mathbf{u} || =1}
 \frac{\widetilde{\mathbf{E}}(\mathbf{u},\mathbf{u})}{1-  \frac{1}{|D|+\rho|Y\setminus D|}
 \left| \int_D \mathbf{u} \right|^2} . \eeq

It is quite interesting to compare \eqnref{crite3} with
\eqnref{crite2}. If $\rho <1$, then
 \beq
 \min_{N_j} \max_{\mathbf{u} \in N_j, || \mathbf{u} || =1}
 \frac{\widetilde{\mathbf{E}}(\mathbf{u},\mathbf{u})}{1 - |\int_D
 \mathbf{u}|^2} < \min_{N_j}\max_{\mathbf{u}\in N_j, || \mathbf{u} || =1}
 \frac{\widetilde{\mathbf{E}}(\mathbf{u},\mathbf{u})}{1-  \frac{1}{|D|+\rho|Y\setminus D|}
 \left| \int_D \mathbf{u} \right|^2},
 \eeq
which shows that smaller the density $\rho$ is, wider the band-gap is, provided that
\eqnref{crite1} is fulfilled. This phenomenon was reported by
Economou and Sigalas \cite{ES94} who observed that periodic
elastic composites whose matrix has lower density and higher shear
modulus compared to those of inclusions yield better open gaps,
and the analysis of this paper agrees with it.

\section{Conclusion}
In this paper, we have reduced band structure calculations for
phononic crystals to the problem of finding the characteristic
values of a family of meromorphic integral operators. We have also
provided complete asymptotic expansions of these characteristic
values as the Lam\'e parameter $\mu$ goes to infinity, established
a connection between the band gap structure and the Dirichlet
eigenvalue problem for the Lam\'e operator, and given a criterion
for gap opening as $\mu$ becomes large. The leading-order terms in
the expansions of the characteristic values are explicitly
computed. An asymptotic analysis for the band-gap structure in
three-dimensions can be provided with only minor modifications of
the techniques presented here. Our results in this paper open the
road to numerous numerical and analytical investigations on
phononic crystals and could, in particular, be used for systematic
optimal design of phononic structures.

\end{document}